\documentclass{amsart}
\usepackage{amssymb}
\usepackage{upref}

\usepackage[all,cmtip]{xy}

\usepackage[lite,initials,shortalphabetic]{amsrefs}

\providecommand{\texorpdfstring}[2]{#1}

\makeatletter
  \renewcommand{\p@enumi}{\thesubsection}
\makeatother

\newenvironment{resumeenumerate}[1]
{\begin{enumerate}
 \setcounter{enumi}{#1}
 \addtocounter{enumi}{-1}
}
{\end{enumerate}
}

\newenvironment{lettered}
{\begin{list}{\thelettercounter)}
 {\usecounter{lettercounter}\def\makelabel##1{\hss\llap{##1}}}
}
{\end{list}
}
\newcounter{lettercounter}
\renewcommand{\thelettercounter}{\alph{lettercounter}}

\newenvironment{resumelettered}[1]
{\begin{lettered}
 \setcounter{lettercounter}{#1}
 \addtocounter{lettercounter}{-1}
}
{\end{lettered}
}

\theoremstyle{plain}

\makeatletter
\newcommand{\emsection}[1]{%
  \par
  \addpenalty\@secpenalty
  \vskip 6 pt plus 9 pt
  \emph{#1.}\nobreak\enspace\ignorespaces
}
\makeatother

\newcommand{\intro}{%
  \goodbreak
  \vskip 6 pt plus 9 pt
}

\newcommand{\boldcentered}[1]{%
  \bigbreak\begin{center}\textbf{#1}\end{center}%
}

\numberwithin{equation}{subsection}

\newcommand{\Comma}{\rlap{\enspace ,}}

\newcommand{\cat}[1]{\boldsymbol{#1}}

\newcommand{\bs}{\boldsymbol}

\newcommand{\RelCat}{\mathbf{RelCat}}

\newcommand{\RelkCat}{\mathbf{Rel}^{k}\mathbf{Cat}}
\newcommand{\RelzCat}{\mathbf{Rel}^{0}\mathbf{Cat}}

\newcommand{\Le}{\mathrm{L}}

\newcommand{\Cat}{\mathbf{Cat}}
\newcommand{\Cathat}{\widehat{\mathbf{Cat}}}

\newcommand{\simp}{\mathrm{s}}

\DeclareMathOperator{\N}{N}

\newcommand{\commacat}[2]{\mathchoice
  {(#1\,\mathord\downarrow\,#2)}
  {(#1\,\mathord\downarrow\,#2)}
  {(#1\mathord\downarrow#2)}
  {(#1\mathord\downarrow#2)}
}
\newcommand{\subcommacat}[3]{\mathchoice
  {(#1\,\mathord\downarrow_{#2}\,#3)}
  {(#1\,\mathord\downarrow_{#2}\,#3)}
  {(#1\mathord\downarrow_{#2}#3)}
  {(#1\mathord\downarrow_{#2}#3)}
}

\DeclareMathOperator{\Ho}{Ho}

\DeclareMathOperator{\Gr}{\boldsymbol{Gr}}
\DeclareMathOperator{\und}{und}
\DeclareMathOperator{\holim}{holim}
\DeclareMathOperator{\wholim}{wholim}

\newcommand{\iso}{\approx}

\newcommand{\op}{^{\mathrm{op}}}

\newcommand{\spacedcdots}{{\cdot\;\cdot\;\cdot}}

\newcommand{\zigzag}[5]{#1\negmedspace: #2
  \mathchoice{\longrightarrow}{\to}{\to}{\to}#3
  \mathchoice{\longleftarrow}{\gets}{\gets}{\gets}#4:\negmedspace #5}

\newcommand{\pullback}[3]{#1{\mathord{\times}_{\!#2}}#3}

\newcommand{\abs}[1]{\left |#1\right|}

\begin{document}

\title[Quillen Theorems $B_{n}$ for homotopy pullbacks of
  $(\infty,k)$-categories]{Quillen Theorems $B_{n}$ for homotopy
  pullbacks\\ of $(\infty,k)$-categories}

\author{C. Barwick}
\address{Department of Mathematics, Massachusetts Institute of
  Technology, Cambridge, MA 02139}
\email{clarkbar@math.mit.edu}

\author{D.M. Kan}
\address{Department of Mathematics, Massachusetts Institute of
  Technology, Cambridge, MA 02139}

%\date{January 20, 2011}
\date{\today}

\begin{abstract}
  We extend the Quillen Theorem $B_{n}$ for \emph{homotopy fibers} of
  Dwyer et~al.\ to similar results for \emph{homotopy pullbacks} and
  note that these results imply similar results for zigzags in the
  categories of \emph{relative categories} and \emph{$k$-relative
    categories}, not only with respect to their \emph{Reedy}
  structures but also their \emph{Rezk} structure, which turns them
  into models for the theories of \emph{$(\infty,1)$- and
    $(\infty,k)$-categories} respectively.

  Our main tool for proving this are the \emph{sharp maps} of Hopkins
  and Rezk which, because of their fibration-like properties, we
  prefer to call \emph{fibrillations}.
\end{abstract}

\maketitle

%--------------------------------------------------------------------
%--------------------------------------------------------------------
% The introduction will be section 0 (for no apparent reason...)
\setcounter{section}{-1}
\section{Introduction}
\label{sec:Intro}

%--------------------------------------------------------------------
\subsection{Background}
\label{sec:Bkgrnd}

\begin{enumerate}
\item \label{Bkgrndi} In \cite{Q}*{\S1} Quillen proved his Theorem~B
  which states that, for a functor $f\colon \cat X \to \cat Z$ and an
  object $Z \in \cat Z$, the rather simple \emph{over category}
  $\commacat{f}{Z}$ is a \emph{homotopy fiber} of $f$ if $f$ has a
  certain \emph{property $B_{1}$}.
\item \label{Bkgrndii} This was generalized in \cite{DKS}*{\S6} where
  it was shown that increasingly weaker \emph{properties} $B_{n}$
  ($n>1$) allowed for increasingly less simple description of these
  homotopy fibers as \emph{$n$-arrow overcategories}
  $\subcommacat{f}{n}{Z}$.
\item \label{Bkgrndiii} It was also noted that a sufficient condition
  for a functor $f\colon \cat X \to \cat Z$ to have property $B_{n}$
  was that the category $\cat Z$ had a certain property $C_{n}$.
\end{enumerate}

%--------------------------------------------------------------------
\subsection{The current paper}
\label{sec:CurPap}

Our main results in this paper are the following.
\begin{enumerate}
\item \label{CurPapi} We show that for a zigzag $\zigzag{f}{\cat
    X}{\cat Z}{\cat Y}{g}$ between categories in which $f$ has
  property $B_{n}$ \eqref{Bkgrndii} (and in particular if $\cat Z$ has
  property $C_{n}$) \eqref{Bkgrndiii}, its \emph{homotopy pullback}
  admits a description similar to the one mentioned in
  \ref{sec:Bkgrnd} namely as a \emph{$n$-arrow pullback category
    $\subcommacat{f\cat X}{n}{g\cat Y}$}.  Moreover its
  \emph{pullback} comes with a monomorphism into the homotopy pullback
  and hence is itself a homotopy pullback if the monomorphism is a
  weak equivalence.
\item \label{CurPapii} We then deduce from this similar results for
  zigzags in the categories of \emph{relative categories} and
  \emph{$k$-relative categories} ($k>1$) and do this not only with
  respect to their \emph{Reedy structure}, but also with respect to
  their \emph{Rezk structure} which turns them into models for the
  theories of \emph{$(\infty,1)$- and $(\infty,k)$-categories}
  respectively.
\item \label{CurPapiii} We also note that a sufficient condition for a
  category, a relative category or a $k$-relative category to have
  \emph{property $C_{3}$} is that it admits what we will call a
  \emph{strict $3$-arrow calculus}.
\item \label{CurPapiv} Our main tool for proving all this consists of
  the \emph{sharp maps} of Hopkins and Rezk \cite{R2} which, because
  of their fibration-like properties we prefer to call
  \emph{fibrillations}.  They are the dual of Hopkins' \emph{flat
    maps} which have similar cofibration-like properties and which we
  therefore call \emph{cofibrillations}.  These cofibrillations do not
  play any role in the current paper, except for a surprise appearance
  in \ref{def:strcalculus}(iii)$'$.
\end{enumerate}

%--------------------------------------------------------------------
\subsection{The genesis of the current paper}
\label{sec:genesis}

The original version of this paper consisted of only the results
mentioned in \ref{CurPapi} and the corresponding part of
\ref{CurPapiii}.  That was exactly what we needed in \cite{BK3} where
it enabled us to reduce the proof that certain pullbacks were
homotopy pullbacks to a rather straightforward calculation.  However
our proof of these results was rather ad hoc and not very
satisfactory.

Fortunately two things happened.
\begin{enumerate}
\item \label{genesisi} We discovered a manuscript of Charles Rezk
  \cite{R2} in which he studied the \emph{sharp maps} of Mike Hopkins.
  These maps seemed to be exactly what we needed.  Just like
  fibrations, they could be used in a \emph{right proper} model
  category to obtain pullbacks which were homotopy pullbacks (which by
  the way made us call them by the more suggestive name of
  \emph{fibrillations}).
\item \label{genesisii} Moreover when subsequently we took a closer
  look at the lemma of Quillen \cite{Q}*{\S1} which started it all and
  which he used to prove his Theorem~B, we noted that this lemma was
  essentially just an elegant way of constructing fibrillations, some
  of which were exactly the ones we needed.
\end{enumerate}

Combining all this with some simple properties of fibrillations we
then not only obtained a much better proof of the above mentioned
result \ref{CurPapi}, but also realized that with relatively little
effort this proof could be extended to a proof of similar results for
\emph{relative categories} and the \emph{$k$-relative categories} of
\cite{BK2}, which as result the current manuscript.

%--------------------------------------------------------------------
\subsection{Acknowledgements}
\label{sec:Ack}

We would like to thank Phil Hirschhorn, Mike Hopkins and especially
Dan Dugger for several useful conversations.

%--------------------------------------------------------------------
%--------------------------------------------------------------------
\setcounter{tocdepth}{1}
\tableofcontents

%--------------------------------------------------------------------
%--------------------------------------------------------------------
\section{An overview}
\label{sec:Oview}

The paper consists of three parts, each of which consists of three
sections.

\boldcentered{Part I: Categorical preliminaries}

In sections 2, 3 and 4 we discuss the various models and relative
categories involved and some of the relative functors between them.
%--------------------------------------------------------------------
\subsection{The categories involved}
\label{sec:CatInvlv}

\begin{enumerate}
\item \label{CatInvlvi} The categories $\RelkCat$ ($k \ge 0$) of small
  \emph{$k$-relative categories} of \cite{BK2}.

  For $k=1$ this is just the category $\RelCat$ of small
  \emph{relative categories} of \cite{BK1}.

  For $k=0$ this is the category of small \emph{maximal} relative
  categories (i.e.\ in which all maps are weak equivalences).  As
  $\RelzCat$ thus is isomorphic, although not equal, to the category
  $\Cat$ of small categories, we \emph{will often denote $\RelzCat$ by
  $\Cathat$}.
\item \label{CatInvoviii} The categories $\simp^{k}\Cathat$ and
  $\simp^{k}\cat S$ ($k\ge 0$) which are respectively the categories
  of \emph{$k$-simplicial diagrams} in $\Cathat$ and the categories of
  small \emph{$k$-simplicial spaces}, i.e.\ $(k+1)$-simplicial sets.
\end{enumerate}

%--------------------------------------------------------------------
\subsection{The functors involved}
\label{sec:FuncInvov}

\begin{enumerate}
\item \label{FuncInvovi} The \emph{$k$-simplicial nerve functors}
  $\simp^{k}\N\colon \RelkCat \to \simp^{k}\cat S$ ($k\ge 0$) of
  \cite{BK2}, which for $k=0$ is just he \emph{classical nerve
    functor $\N$}.
\item \label{FuncInfovii} The functors $\N_{*}\colon \simp^{k}\Cathat
  \to \simp^{k}\cat S$ ($k\ge 0$) induced by $\N$.
\item The (unique) functors $w_{*}\colon \RelkCat \to
  \simp^{k}\Cathat$ ($k\ge 0$) such that $\N_{*}w_{*} = \simp^{k}\N$.
\end{enumerate}

%--------------------------------------------------------------------
\subsection{The model and relative structures}
\label{sec:ModRlSt}

\begin{enumerate}
\item \label{ModRlSti} We endow $\cat S = \simp^{0}\cat S$ with the
  classical model structure and $\Cathat$ with the Quillen equivalent
  Thomason structure \cite{T2}.
\item \label{ModRlStii} We endow $\simp^{k}\cat S$ and
  $\simp^{k}\Cathat$ ($k\ge 1$) with the resulting Reedy model structure.
\item \label{ModRlStiii} We endow $\RelCat$ and $\RelkCat$ ($k>1$)
  respectively with the Quillen equivalent Reedy \emph{model}
  structure which \cite{BK1} is lifted from the Reedy structure on
  $\simp\cat S$ and the homotopy equivalent \emph{relative} Reedy
  structure which \cite{BK2} is lifted from the Reedy structure on
  $\simp^{k}\cat S$.
\item \label{ModRlStiv} In the application of our main result we will
  consider the categories $\RelkCat$ and $\simp^{k}\Cathat$ ($k\ge 1$)
  as models for the theory of $(\infty,k)$-categories by endowing them
  with a \emph{Rezk} structure which has more weak equivalences than
  the Reedy structure and we will denote those categories so endowed
  by $\Le\RelkCat$ and $\Le\simp^{k}\Cathat$.
\end{enumerate}

%--------------------------------------------------------------------
\subsection{Properties of the functors involved}
\label{sec:FncPrp}

The functors considered in \ref{sec:FuncInvov} above all have the same
nice properties as the classical nerve functors (and we will therefore
refer to them as \emph{abstract nerve functors}).

In particular they all are \emph{relative functors} which are
\emph{homotopy equivalences}.

%--------------------------------------------------------------------
%--------------------------------------------------------------------
\boldcentered{Part II: Homotopy pullbacks and potential homotopy
  pullbacks}

\subsection{}
\label{sec:empty}

In section \ref{sec:HmtpyPlbk} we will define homotopy pullbacks of
zigzags in a model category in a more general fashion than is usually
done, but that enables us to extend the definition to certain
saturated relative categories.
\begin{enumerate}
\item \label{emptyi} In a model category we define a homotopy pullback
  of a zigzag as \emph{any} object which is weakly equivalent to its
  image under a \emph{``homotopically correct'' homotopy limit functor}.
\item \label{emptyii} In a saturated relative category we then define
  a homotopy pullback of a zigzag of a \emph{any} object which is
  weakly equivalent to its image under what we will call a \emph{weak
    homotopy limit functor} which is a functor which has only some of
  the properties of the above (i) homotopy limit functors.
\item \label{emptyiii} Our main result then is a \emph{Global
    equivalence lemma} which states that
  \begin{itemize}
  \item \emph{if $f\colon \cat C \to \cat D$ is a homotopy equivalence
    between saturated relative categories,} then \emph{$\cat C$ admits
  weak homotopy limit functors iff $\cat D$ does, and in that case $f$
  preserves homotopy pullbacks.}
  \end{itemize}
\item \label{emptyiv} This lemma not only
  \begin{itemize}
  \item implies that any two weak homotopy limit functors on the same
    category yield the same notion of homotopy pullbacks of zigzags, and
  \item enables us, in view of the fact that model categories admit
    such weak homotopy limit functors, to prove their existence
    elsewhere.
  \end{itemize}
  but also plays a role in the proof of our main result in section
  \ref{sec:PrfThm}.
\end{enumerate}

%--------------------------------------------------------------------
\subsection{}
\label{sec:onep6}

In section \ref{sec:PotHmPb} we discuss, for a zigzag $\zigzag{f}{\cat
  X}{\cat Z}{\cat Y}{g}$ in $\RelkCat$ or $\simp^{k}\Cathat$ ($k\ge
0$) the \emph{$n$-arrow pullback object} $\subcommacat{f\cat
  X}{n}{g\cat Y}$ which was mentioned in \ref{sec:CurPap} above.

%--------------------------------------------------------------------
\subsection{}
\label{sec:onep7}

In section \ref{sec:PrpBnCn} we recall the notions of \emph{property
  $B_{n}$} and \emph{property $C_{n}$} for $\Cathat$ which were
mentioned in \ref{sec:Bkgrnd} above and then extend them to the
categories $\RelkCat$ and $\simp^{k}\Cathat$ for $k\ge 1$.

%--------------------------------------------------------------------
%--------------------------------------------------------------------
\boldcentered{Part III: The main results and their proof}

\subsection{}
\label{sec:onep8}

In section \ref{sec:MainRlt} we formulate our main results, Theorems
\ref{thm:MpPropBn}, \ref{thm:SuffCn}, \ref{thm:QFibQfib} and
\ref{thm:MpBn}.

We also prove Theorem \ref{thm:SuffCn} and, assuming Theorem
\ref{thm:MpPropBn}, Theorems \ref{thm:QFibQfib} and \ref{thm:MpBn} as
well, but defer the proof of Theorem \ref{thm:MpPropBn} itself until
section \ref{sec:PrfThm}.

\begin{enumerate}
\item \label{onep8i} Theorem \ref{thm:MpPropBn} states that, for a
  zigzag $\zigzag{f}{\cat X}{\cat Z}{\cat Y}{g}$ in $\RelkCat$ or
  $\simp^{k}\Cat$ ($k\ge 0$), for which the map $f\colon \cat X \to
  \cat Z$ has property $B_{n}$ or the object $\cat Z$ has property
  $C_{n}$ \eqref{sec:onep7}, the \emph{$n$-arrow pullback object}
  $\subcommacat{f\cat X}{n}{g\cat Y}$ \eqref{sec:onep6} is actually a
  \emph{homotopy pullback} \eqref{sec:empty} of that zigzag.
\item \label{onep8ii} Theorem \ref{thm:MpPropBn} states that a
  sufficient condition on an object $\cat Z \in \RelkCat$ ($k\ge 0$)
  in order that it has property $C_{3}$ is that $\cat Z$ admits what
  we will call a \emph{strict $3$-arrow calculus}.
\item \label{onep8iii} If $k=1$, then this (ii) is in particular the
  case if $\cat Z$ is a \emph{partial model category}, i.e.\ a
  relative category which has the \emph{two out of six} property and
  admits a \emph{$3$-arrow calculus}.

  These partial model categories have the property that \cite{BK3}, if
  $\RelCat$ is endowed with the \emph{Rezk} structure
  \eqref{ModRlStiv}, then
  \begin{itemize}
  \item every partial model category is Reedy equivalent to a fibrant
    object in $\Le\RelCat$ \eqref{ModRlStiv}, and
  \item every fibrant object in $\Le\RelCat$ is Reedy equivalent to a
    partial model category.
  \end{itemize}
  Consequently one has the following result for
  $(\infty,1)$-categories.
\item \label{onep8iv} Theorem \ref{thm:QFibQfib} that states that, for
  every zigzag $\zigzag{f}{\cat X}{\cat Z}{\cat Y}{g}$ in $\RelCat$
  between partial model categories, the \emph{$3$-arrow pullback
    object} $\subcommacat{f\cat X}{3}{g\cat Y}$ is a \emph{homotopy
    pullback} of this zigzag not only in $\RelCat$, but also in
  $\Le\RelCat$, i.e.\ in $(\infty,1)$-categories.
\item \label{onep8v} Finally in Theorem \ref{thm:MpBn} we state a
  similar result for $\Le\RelkCat$ ($k>1$), i.e.\ for
  $(\infty,k)$-categories, which however is much weaker than that
  Theorem \ref{thm:QFibQfib}, as we have no model structure on
  $\RelkCat$ for $k>1$ nor an analog for partial model categories.
\end{enumerate}

%--------------------------------------------------------------------
\subsection{}
\label{sec:onep9}

In section \ref{sec:HopRzF} we review the \emph{fibrillations} of
Hopkins and Rezk which can be defined in any \emph{relative category
  with pullbacks} and describe the following three lemmas which we
need in section \ref{sec:PrfThm} to prove Theorem \ref{thm:MpPropBn}.

\begin{enumerate}
\item \label{onep9i} A \emph{Quillen fibrillation lemma} which is a
  reformulation and a slight strengthening of the lemma that Quillen
  used to prove his Theorem~B, and which enables us to obtain the
  needed fibrillations in $\Cathat$.

  If, for an object $\cat D \in \Cathat$ and a not necessarily
  relative functor $F\colon \cat D \to \Cathat$, $\Gr F \in \Cathat$
  denotes its Grothendieck construction and $\pi\colon \Gr F \to \cat
  D \in \Cathat$ denotes the associated projection functor, then this
  lemma states that
  \begin{itemize}
  \item the projection functor $\pi\colon \Gr F \to \cat D \in
    \Cathat$ is a \emph{fibrillation} iff the functor $F\colon \cat D
    \to \Cathat$ is a \emph{relative} functor, i.e.\ sends every map
    of $\cat D$ to a weak equivalence in $\Cathat$.
  \end{itemize}
\item \label{onep9ii} A \emph{fibrillation lifting lemma} which
  enables us to obtain from the fibrillation in $\Cathat$ obtained in
  (i) above corresponding fibrillations in the categories
  $\simp^{k}\Cathat$ for $k\ge 1$.
\item \label{onep9iii} A \emph{Hopkins-Rezk fibrillation lemma} which
  enables us to use the fibrillations of (ii) above to obtain the
  desired homotopy pullback in the categories $\simp^{k}\Cat$ ($k\ge
  0$), i.e.\ it states that
  \begin{itemize}
  \item for every pullback square in a right proper model category
    \begin{displaymath}
      \xymatrix{
        {\cat U} \ar[r] \ar[d]
        & {\cat Y} \ar[d]\\
        {\cat V} \ar[r]
        & {\cat Z}
      }
    \end{displaymath}
    in which the map $\cat C \to \cat Z$ is a fibrillation, the object
    $\cat U$ is a homotopy pullback of the zigzag $\cat V \to \cat Z
    \gets \cat Y$.
  \end{itemize}
\end{enumerate}

%--------------------------------------------------------------------
\subsection{}
\label{sec:onep10}

In section \ref{sec:PrfThm} we finally prove Theorem
\ref{thm:MpPropBn}.  We prove this first for the categories
$\simp^{k}\Cathat$ ($k\ge 0$) by using the above \eqref{sec:onep9}
\emph{three lemmas} and then lift these results to the categories
$\RelkCat$ ($k\ge 1$) by means of the \emph{Global equivalence lemma}
which was mentioned in \ref{emptyiii} above.

Actually for $\RelCat$ we did not have to use this fourth lemma as we
could have obtained this also using only the first three.

%--------------------------------------------------------------------
%--------------------------------------------------------------------
%--------------------------------------------------------------------
\part{Categorical Preliminaries}

\section{Relative categories}
\label{sec:RelCat}

%--------------------------------------------------------------------
\subsection{Summary}
\label{sec:RlCtSum}
We start with a brief review of
\begin{enumerate}
\item \emph{relative categories} and \emph{relative functors} between them,
\item \emph{saturated}, \emph{maximal} and \emph{minimal} relative
  categories,
\item a \emph{homotopy relation} between relative functors,
\item \emph{strict $3$-arrow calculi} on relative categories,
\end{enumerate}
and
\begin{resumeenumerate}{5}
\item the associated \emph{simplicial nerve functor} to simplicial
  spaces, i.e.\ bisimplicial sets.
\end{resumeenumerate}

%--------------------------------------------------------------------
%--------------------------------------------------------------------
\boldcentered{Relative Categories}

\subsection{Definition}
\label{def:relcat}

\begin{enumerate}
\item \label{relcati} A \textbf{relative category} is a pair $(\cat C,
  w\cat C)$ (usually denoted by just $\cat C$) consisting of an
  \textbf{underlying category} $\cat C$ (sometimes denoted by
  $\und\cat C$) and a subcategory $w\cat C \subset \cat C$ which is
  \emph{only} required to contain \emph{all} the objects of $\cat C$
  (and hence their identity maps).
\item \label{relcatii} The maps or $w\cat C$ are called \textbf{weak
    equivalences} and two objects of $\cat C$ are called
  \textbf{weakly equivalent} if they can be connected by a finite
  zigzag of weak equivalences.
\item \label{relcatiii} By a \textbf{functor} $f\colon \cat C \to \cat
  D$ between two relative categories we will mean just a functor
  $f\colon \und\cat C \to \und\cat D$ (i), and
\item \label{relcativ} such a functor will be called a
  \textbf{relative functor} if it preserves weak equivalences, i.e.\
  if $f(w\cat C) \subset w\cat D$.
\end{enumerate}

%--------------------------------------------------------------------
\boldcentered{Saturated, maximal and minimal relative categories}

%--------------------------------------------------------------------
\subsection{Definition}
\label{def:saturated}

\begin{enumerate}
\item \label{saturi} A relative category $\cat C$ will be called
  \textbf{saturated} if a map in $\cat C$ is a weak equivalence iff
  its image in $\Ho\cat C$ (i.e.\ the category obtained from $\cat C$
  by formally inverting all the weak equivalences) is an isomorphism.
\end{enumerate}
This is in particular the case if
\begin{resumeenumerate}{2}
\item \label{saturii} $\cat C$ is a \textbf{maximal} relative
  category, i.e.\ \emph{all} maps of $\cat C$ are weak equivalences
\end{resumeenumerate}
or if
\begin{resumeenumerate}{3}
  \item \label{saturiii} $\cat C$ is a \textbf{minimal} relative
    category, i.e.\ the identity maps re the \emph{only} weak
    equivalences.
\end{resumeenumerate}

%--------------------------------------------------------------------
\boldcentered{A homotopy relation between relative functors}

%--------------------------------------------------------------------
\subsection{Definition}
\label{def:HomRel}

Let $\bs 1^{w}$ denote the maximal relative category
(Df.~\ref{saturii}) which  consists of two objects $0$ and $1$ and
one map $0\to 1$ between them.  Then $\bs 1^{w}$ gives rise to the
following notions.
\begin{enumerate}
\item \label{HomReli} Given two relative functors $f,g\colon \cat C
  \to \cat D$ between relative categories, a \textbf{strict homotopy}
  $h\colon f \to g$ between them will be a natural weak equivalence,
  i.e.\ a map
  \begin{displaymath}
    h\colon \cat C\times\bs 1^{w} \longrightarrow \cat D
  \end{displaymath}
  such that $h(c,0) = fc$ and $h(c,1) = gc$ for every object or map $c
  \in \cat C$.
\item \label{HomRelii} Two relative functors $\cat C \to \cat D$ then
  are called \textbf{homotopic} if they can be connected by a finite
  zigzag of strict homotopies, and
\item \label{HomReliii} a functor $f\colon \cat C \to \cat D$ is
  called a \textbf{homotopy equivalence} if $f$ is a relative functor
  and if there exists a relative functor $f'\colon \cat D \to \cat C$
  (called a \textbf{homotopy inverse} of $f$) such that the
  compositions $f'f$ and $ff'$ are homotopic to $1_{\cat C}$ and
  $1_{\cat D}$ respectively.
\end{enumerate}

%--------------------------------------------------------------------
\boldcentered{$3$-arrow calculi}

%--------------------------------------------------------------------
\subsection{Definition}
\label{def:calculus}

A relative category $\cat C$ is said to admit a \textbf{$3$-arrow
  calculus} if there exist subcategories $\cat U$ and $\cat V \subset
w\cat C$ which behave like the categories of the \emph{trivial
  cofibrations} and \emph{trivial fibrations} in a model category in
the sense that
\begin{enumerate}
\item \label{calculusi} for every map $u \in \cat U$, its pushouts in
  $\cat C$ exist and are again in $\cat U$,
\item \label{calculusii} for every map $v \in \cat V$, it's pullbacks
  in $\cat C$ exist and are again in $\cat C$, and
\item \label{calculusiii} the maps $w \in w\cat C$ admit a
  \emph{functorial} factorization $w = vu$ with $u \in \cat U$ and $v
  \in \cat V$.
\end{enumerate}

%--------------------------------------------------------------------
\subsection{Remark}
\label{rem:easier}

It should be noted that the conditions \ref{calculusi} and (ii) are
stronger than the ones in \cite{DK}*{8.2} and \cite{DHKS}*{36.1}.
However we prefer them as they are easier to use and are likely to be
automatically satisfied.

%--------------------------------------------------------------------
\subsection{Remark}
\label{rem:whystrict}

In \cite{DK}*{8.2} and \cite{DHKS}*{36.1} $3$-arrow calculi were
defined in the presence of the \emph{two out of three} and the
\emph{two out of six} properties respectively with the result that
\emph{a $3$-arrow calculus on $(\cat C,w\cat C)$ automatically
  restricted to a $3$-arrow calculus on $(w\cat C, w\cat C)$}.

As we do not want to assume the presence of the two out of
three or the two out of six property, however, we define a notion of
\emph{strict $3$-arrow calculi} as follows.

%--------------------------------------------------------------------
\boldcentered{Strict $3$-arrow calculi}

%--------------------------------------------------------------------
\subsection{Definition}
\label{def:strcalculus}

A \textbf{strict $3$-arrow calculus} on a relative category $(\cat
C,w\cat C)$ will be a $3$-arrow calculus (Df.~\ref{def:calculus})
which restricts to a $3$-arrow calculus on $(w\cat C, w\cat C)$.

Another way of saying this by adding in Df.~\ref{def:calculus} to the
conditions (i) and (ii) the conditions that
\begin{enumerate}
\item [(i)$'$] every pushout of a map in $w\cat C$ along a map in
  $\cat U$ is again in $w\cat C$, and
\item [(ii)$'$] every pullback of a map in $w\cat C$ along a map in
  $\cat V$ is again in $w\cat C$.
\end{enumerate}

With other words the maps in $\cat U$ and $\cat V$ behave like trivial
cofibrations and trivial fibrations in a \emph{proper} model category.

A more compact way of saying all this is that
\begin{enumerate}
\item [(iii)$'$] a strict $3$-arrow calculus on a relative category is
  a functorial factorization of its weak equivalences into a
  \emph{trivial cofibrillation} followed by a \emph{trivial fibrillation}
\end{enumerate}
where fibrillations are as defined in Df.~\ref{def:fibril} below and
cofibrillations are defined dually.

%--------------------------------------------------------------------
\boldcentered{The simplicial nerve functor}

%--------------------------------------------------------------------
\subsection{Definition}
\label{def:smpnv}

Let $\RelCat$ denote the category of (small) \emph{relative
  categories} (Df.~\ref{def:relcat}) and let $\simp\cat S$ denote the
category of \textbf{simplicial spaces}, i.e.\ \emph{bisimplicial}
sets.

Furthermore, for every integer $p \ge 0$, let $\cat p$ denote the
$p$-arrow category
\begin{displaymath}
  0 \longrightarrow \cdots \longrightarrow p
\end{displaymath}
and let $\cat p^{v}$ and $\cat p^{w}$ denote respectively the
\emph{minimal} and \emph{maximal} relative categories
(Df.~\ref{ModRlStiii} and (ii)) which have $\cat p$ as their
underlying category.

The \textbf{simplicial nerve functor} then is the functor
\begin{displaymath}
  \simp\N\colon \RelCat \longrightarrow \simp\cat S
\end{displaymath}
which sends each object $\cat C \in \RelCat$ to the bisimplicial set
which in bidimension $(p,q)$ consists of the maps
\begin{displaymath}
  \cat p^{v}\times \cat q^{w} \longrightarrow \cat C \in \RelCat
\end{displaymath}

%--------------------------------------------------------------------
%--------------------------------------------------------------------
\section{The categories and functors involved}
\label{sec:CtFncInv}

%--------------------------------------------------------------------
\subsection{Summary}
\label{sec:CtFncSumm}

We now describe the categories and functors which we need in the paper
and note that
\begin{enumerate}
\item \label{CtFncSummi} these categories come with obvious notions of
  \emph{strict homotopies} and \emph{homotopy equivalences} (as in
  Df.~\ref{rem:easier}),
\item \label{CtFncSummii} these functors \emph{preserve} these strict
  homotopies and homotopy equivalences, and
\item \label{CtFncSummiii} these functors all have \emph{left
    adjoints} which are \emph{left inverses}.
\end{enumerate}

In section \ref{sec:AbsNv} we then endow these categories with model
or relative structures which
\begin{resumeenumerate}{4}
\item \label{CtFncSummiv} are compatible with these homotopy
  equivalences in the sense that \emph{every homotopy equivalence is a
  weak equivalence}, and
\item \label{CtFncSummv} turn the functors between these categories
  into relative functors which are \emph{homotopy equivalences}.
\end{resumeenumerate}

%--------------------------------------------------------------------
%--------------------------------------------------------------------
\boldcentered{$k$-relative categories}

%--------------------------------------------------------------------
\subsection{Definition}
\label{def:krel}

We extend the definition of $k$-relative categories for $k \ge 1$ of
\cite{BK2}*{2.3} to include also the case $k=0$ as follows.

A \textbf{$k$-relative category} ($k>0$) will be a $(k+2)$-tuple
\begin{displaymath}
  \cat C = (a\cat C, v_{1}\cat C, \ldots, v_{k}\cat C, w\cat C)
\end{displaymath}
consisting of an \textbf{ambient category} $a\cat C$ and subcategories
\begin{displaymath}
  v_{1}\cat C, \ldots, v_{k}\cat C \text{ and } w\cat C
  \subset a\cat C
\end{displaymath}
that each contain all the objects of $a\cat C$ and form a
commutative diagram with $2k$ arrows of the form
\begin{displaymath}
  \xymatrix@C=0.5em{
    & {w\cat C} \ar[dl] \ar[dr]\\
    {v_{1}\cat C} \ar[dr]
    & {\spacedcdots}
    & {v_{k}\cat C} \ar[dl]\\
    & {a\cat C}
  }
\end{displaymath}
and where $a\cat C$ is subject to the conditions that
\begin{enumerate}
\item \label{kreli} every map in $a\cat C$ is a finite composition of
  maps in the $v_{i}\cat C$ ($1 \le i \le k$), and
\item \label{krelii} every relation in $a\cat C$ is the consequence of
  the commutativity relations in $v_i\cat C$ ($1\leq i\leq n$) and the commutativity in $a\cat C$ of diagrams of the form
  \begin{displaymath}
    \xymatrix{
      {\cdot} \ar[r]^{x_{1}} \ar[d]_{y_{1}}
      & {\cdot} \ar[d]^{y_{2}}\\
      {\cdot} \ar[r]_{x_{2}}
      & {\cdot}
    }
  \end{displaymath}
  in which $x_{1}, x_{2} \in v_{i}\cat C$ and $y_{1}, y_{2} \in
  v_{j}\cat C$ and $1\leq i<j\leq n$.
\end{enumerate}

The maps of $w\cat C$ are called \textbf{weak equivalences} and a
\textbf{relative functor} $f\colon \cat C \to \cat D$ between two
$k$-relative categories will be a functor $f\colon a\cat C \to a\cat
D$ such that
\begin{displaymath}
  f(w\cat C) \subset w\cat D
  \qquad\text{and}\qquad
  f(v_{i}\cat C) \subset v_{i}\cat D
  \qquad \text{for all $i \le i < k$.}
\end{displaymath}

%--------------------------------------------------------------------
\subsection{Remark}
\label{rem:onerel}

If $\cat C$ is a $1$-relative category, then $a\cat C = v_{1}\cat C$
and hence
\begin{itemize}
\item \emph{$1$-relative categories are the same as relative
    categories} (Df.~\ref{def:relcat}).
\end{itemize}

%--------------------------------------------------------------------
\subsection{Remark}
\label{rem:zerorel}

If $\cat C$ is a $0$-relative category, then $a\cat C = w\cat C$ and
hence
\begin{itemize}
\item $0$-relative categories are the same as maximal relative
  categories (Df.~\ref{saturii}).
\end{itemize}

%--------------------------------------------------------------------
%--------------------------------------------------------------------
\boldcentered{Strict $3$-arrow calculi and $k$-relative categories}

%--------------------------------------------------------------------
\subsection{Definition}
\label{def:Str3ar}

A \textbf{strict $3$-arrow calculus} on a $k$-relative category $\cat
C$ is a $3$-arrow calculus (Df.~\ref{def:calculus}) on $(a\cat C,w\cat
C)$ which restricts to a $3$-arrow calculus on
\begin{displaymath}
  (w\cat C,w\cat C)\qquad\text{and}\qquad (v_{i}\cat C,w\cat C)
  \qquad \text{for all $i \le i \le k$.}
\end{displaymath}

%--------------------------------------------------------------------
%--------------------------------------------------------------------
\boldcentered{The $k$-simplicial nerve functor}

%--------------------------------------------------------------------
\subsection{Definition}
\label{def:kSmpNv}

Let $\RelkCat$ ($k\ge 0$) denote the category of (small) $k$-relative
categories (Df.~\ref{def:krel}) and let $\simp^{k}\cat S$ denote the
category of (small) \textbf{$k$-simplicial spaces}, i.e.\
$(k+1)$-simplicial sets.

Furthermore, for every integer $p \ge 0$ let $\abs{\cat p} \subset
\cat p$ (Df.~\ref{def:smpnv}) denote the subcategory which consists of
the objects and their identity maps only and let
\begin{displaymath}
  \text{$\cat p^{w}$ and
  $\cat p^{v_{i}}$ ($1 \le i \le k$) $\in \RelkCat$}
\end{displaymath}
be such that
\begin{gather*}
  a\cat p^{w} = v_{i}\cat p^{w} = w\cat p^{w} = \cat p
  \qquad \text{($1 \le i \le k$), and}\\
  a\cat p^{v_{i}} = v_{i}\cat p^{v_{i}} = \cat p
  \quad\text{and}\quad
  v_{j}\cat p^{v_{i}} = w\cat p^{v_{i}} = \abs{\cat p}
  \quad\text{for $j \ne i$.}
\end{gather*}

The \textbf{$k$-simplicial nerve functor} then is the functor
\begin{displaymath}
  \simp^{k}\N\colon \RelkCat \longrightarrow \simp^{k}\cat S
\end{displaymath}
which each object $\cat C \in \RelkCat$ to the $(k+1)$-simplicial set
which in multidimension $(p_{k}, \ldots, p_{1}, q)$ consists of the
maps
\begin{displaymath}
  \cat p_{k}^{v_{k}} \times \cdots \times \cat p_{1}^{v_{1}}\times
  \cat q^{w} \longrightarrow \cat C \in \RelkCat
\end{displaymath}

%--------------------------------------------------------------------
\subsection{Remark}
\label{rem:agree}

If $k=1$, then (Rk.~\ref{rem:onerel}) $\RelkCat = \RelCat$
(Df.~\ref{def:smpnv}) and hence in that case Df.~\ref{def:kSmpNv}
agrees with Df.~\ref{def:smpnv}.

%--------------------------------------------------------------------
\subsection{Notation}
\label{not:kzero}

If $k=0$, then
\begin{enumerate}
\item \label{kzeroi} the category $\RelzCat$ is, in view of
  Rk.~\ref{rem:zerorel}, isomorphic, although not identical, with the
  category $\Cat$ of (small) categories and we will therefore
  \begin{itemize}
  \item \emph{often denote $\RelzCat$ by $\Cathat$.}
  \end{itemize}
\item \label{kzeroii} $\simp^{0}\cat S = \cat S$, the category of
  (small) simplicial sets, and
\item \label{kzeroiii} the functor $\simp^{0}\N\colon \RelzCat \to
  \simp^{0}\cat S$ is just the \textbf{classical nerve functor} which
  we will therefore denote by
  \begin{displaymath}
    \N\colon \Cathat \longrightarrow \cat S
  \end{displaymath}
\end{enumerate}

Closely related to the $k$-simplicial nerve functor are the following
two functors.

%--------------------------------------------------------------------
%--------------------------------------------------------------------
\boldcentered{The levelwise nerve functor}

%--------------------------------------------------------------------
\subsection{Definition}
\label{def:lvlnrv}

Let $k\ge 0$.  Then the \textbf{levelwise nerve functor}
\begin{displaymath}
  \N_{*}\colon \simp^{k}\Cathat \longrightarrow \simp^{k}\cat S
\end{displaymath}
will be the functor in which $\simp^{k}\Cathat$ denotes the category
of $k$-simplicial diagrams in $\Cathat$ and $\N_{*}$ is dimensionwise
application of the classical nerve functor $\N$
(Nt.~\ref{rem:onerel}).

%--------------------------------------------------------------------
%--------------------------------------------------------------------
\boldcentered{The higher equivalence functor}

%--------------------------------------------------------------------
\subsection{Definition}
\label{def:HighEq}

For every integer $k \ge 0$ we denote by
\begin{displaymath}
  w_{*}\colon \RelkCat \longrightarrow \simp^{k}\Cathat
\end{displaymath}
the \textbf{higher equivalence functor} which, for every $k$-fold
dimension $p_{*} = (p_{k},\ldots,p_{1})$ sends each object $\cat C \in
\RelkCat$ to the object $w_{p_{*}}\cat C \in \Cathat$ which has as its
objects the maps
\begin{displaymath}
  \cat p_{k}^{v_{k}} \times \cdots \times \cat p_{1}^{v_{1}}
  \longrightarrow \cat C \in \RelkCat
\end{displaymath}
and as maps the natural weak equivalences between them.

\intro The above functors have the following properties.
%--------------------------------------------------------------------
\subsection{Proposition}
\label{prop:RlCtsmpS}

\emph{The following diagram commutes}
\begin{displaymath}
  \vcenter{
    \xymatrix{
      {\RelkCat} \ar[rr]_-{\simp^{k}\N} \ar[dr]_{w_{*}}
      && {\simp^{k}\cat S}\\
      & {\simp^{k}\Cathat} \ar[ur]_{\N_{x}}
    }
  }
  \qquad\qquad
  k\ge 0
\end{displaymath}

\begin{proof}
  This follows directly from the definitions.
\end{proof}

%--------------------------------------------------------------------
\subsection{Proposition}
\label{prop:SmpCtAdj}

\begin{em}
  For every integer $k \ge 0$ the functors
  \begin{enumerate}
  \item $\N_{*}\colon \simp^{k}\Cathat \to \simp^{k}\cat S$,
  \item $\simp^{k}\N\colon \RelkCat \to \simp^{k}\cat S$, and
  \item $w_{*}\colon \RelkCat \to \simp^{k}\Cathat$
  \end{enumerate}
  have a left adjoint which is a left inverse.
\end{em}

\begin{proof}
  \leavevmode
  \begin{enumerate}
  \item It clearly suffices to prove this for the functor $\N$.  But
    this was already done long ago in \cite{GZ}*{II,~4.7}.
  \item This was shown in \cite{BK2}*{4.3}.
  \item This follows from Pr.~\ref{prop:RlCtsmpS} above and the
    following result.  \qedhere
  \end{enumerate}
\end{proof}

%--------------------------------------------------------------------
\subsection{Proposition}
\label{prop:LfAdjInvComp}

If, for two functors $g_{1}\colon \cat A \to \cat B$ and $g_{2}\colon
\cat B \to \cat C$, both $g_{2}$ and $g_{2}g_{1}$ have left adjoints
which are left inverses, then so does $g_{1}$.

\begin{proof}
  Consider the diagram
%   \begin{displaymath}
%     \xymatrix{
%       {\cat A} \ar@<-0.25ex>[r]_{g_{1}} \ar@<-0.75ex>[r];[]_{f_{1}}
%       & {\cat B} \ar@<-0.25ex>[r]_{g_{2}} \ar@<-0.75ex>[r];[]_{f_{2}}
%       & {\cat C} \ar@(u,u)[ll]^{f_{12}}
%     }
%   \end{displaymath}
  \begin{displaymath}
    \xymatrix{
      {\cat A} \ar@<-0.25ex>[r]_{g_{1}} \ar@<-0.75ex>[r];[]_{f_{1}}
      & {\cat B} \ar@<-0.25ex>[r]_{g_{2}} \ar@<-0.75ex>[r];[]_{f_{2}}
      & {\cat C} \ar `u[ll] `[ll]_{f_{12}} [ll]
    }
  \end{displaymath}
  in which $f_{2}$ and $f_{12}$ are left adjoints of $g_{2}$ and
  $g_{2}g_{1}$ and
  \begin{displaymath}
    f_{2}g_{2} = 1,\qquad
    f_{12}g_{2}g_{1} = 1
    \qquad\text{and}\qquad
    f_{1} = f_{12}g_{2}.
  \end{displaymath}
  Then $f_{1}g_{1} = f_{12}g_{2}g_{1} = 1$ and it thus remains to show
  that $f_{1}$ is a left adjoint of $g_{1}$ and this one does by
  noting that every pair of objects $A \in \cat A$ and $B \in \cat B$
  gives rise to a natural sequence of identity maps and isomorphisms
%   \begin{displaymath}
%     \cat A(f_{1}B,A)
%     = \cat A(f_{12}g_{2}B,A)
%     \iso \cat C(g_{2}B, g_{2}g_{1}A)
%     \iso \cat B(f_{2}g_{2}B,g_{1}A)
%     = \cat B(B,g_{1}A)
%   \end{displaymath}
  \begin{displaymath}
    \vcenter{
    \xymatrix@C=1em@R=1.5ex{
      {\cat A(f_{1}B,A)} \ar@{=}[d]
      && {\cat B(B,g_{1}A)} \ar@{=}[d]\\
      {\cat A(f_{12}g_{2}B,A)} \ar@{}[r]|{\iso}
      & {\cat C(g_{2}B, g_{2}g_{1}A)} \ar@{}[r]|{\iso}
      & {\cat B(f_{2}g_{2}B,g_{1}A)}
    }
    }  \qedhere
  \end{displaymath}
\end{proof}

\intro
It remains to deal with the results that were promised in 
\ref{CtFncSummi} and (ii).
%--------------------------------------------------------------------
%--------------------------------------------------------------------
\boldcentered{Homotopy relations on $\RelkCat$, $\simp^{k}\cat S$ and
    $\simp^{k}\Cat$}

%--------------------------------------------------------------------
\subsection{Definition}
\label{def:StrHmtp}

One can define \textbf{strict homotopies} and \textbf{homotopy
  equivalences} in the categories $\RelkCat$, $\simp^{k}\cat S$ and
$\simp^{k}\Cathat$ as in Df.~\ref{def:HomRel} by choosing an ``obvious
analog'' of $\bs 1^{w}$ as follows.
\begin{enumerate}
\item \label{StrHmtpi} For $\RelkCat$ this will be the object $\bs
  1^{w} \in \RelkCat$ (Df.~\ref{def:kSmpNv}),
\item \label{StrHmtpii} for $\simp^{k}\cat S$ this will be the
  standard multisimplex $\Delta[0,\ldots,0,1] \in \simp^{k}\cat S$,
  and
\item \label{StrHmtpiii} for $\simp^{k}\Cathat$ this will be the
  object $\bs 1^{w}_{x} \in \simp^{k}\Cathat$ such that, for every
  $k$-fold dimension $p_{*} = (p_{k}, \ldots, p_{1})$, $\bs
  1^{w}_{p_{*}} = \bs 1^{w} \in \Cathat$.
\end{enumerate}

%--------------------------------------------------------------------
\subsection{Proposition}
\label{prop:PrsStHmt}

\begin{em}
  The functors $\simp^{k}\N$, $w_{*}$ and $\N_{*}$
  \begin{enumerate}
  \item send strictly homotopic maps to strictly homotopic maps
  \end{enumerate}
  and hence
  \begin{resumeenumerate}{2}
    \item send homotopy equivalences to homotopy equivalences.
  \end{resumeenumerate}
\end{em}

\begin{proof}
  \leavevmode
  \begin{enumerate}
  \item The classical nerve functor $\N\colon \Cathat \to \cat S$ has
    the property that (\cite{La}*{\S2} and \cite{Le}), for every two
    maps $f,g\colon \cat A \to \cat B \in \Cathat$
    \begin{itemize}
    \item $\N$ sends every strict homotopy between $f$ and $g$ to a
      strict homotopy between $\N f$ and $\N g$, and
    \item conversely, every strict homotopy in $\cat S$ between $\N f$
      and $\N g$ is obtained in this fashion.
    \end{itemize}
  \item This readily implies that $\N_{*}$ has the same properties.
  \item The proof for $\simp^{k}\N$ is just a higher dimensional
    version of the proof for $k=1$ in \cite{BK1}*{7.5}.
  \item A proof for $w_{*}$ then is obtained by combining
    Pr.~\ref{prop:RlCtsmpS} with (ii) and (iii) above.\qedhere
  \end{enumerate}
\end{proof}

%--------------------------------------------------------------------
%--------------------------------------------------------------------
\section{Abstract nerve functors}
\label{sec:AbsNv}

%--------------------------------------------------------------------
\subsection{Summary}
\label{sec:AbsNvSm}

We now endow the categories $\RelkCat$, $\simp^{k}\cat S$ and
$\simp^{k}\Cathat$ ($k\ge 0$) with a \emph{relative} or \emph{model}
structure which will turn the functors $\simp^{k}\N$, $w_{*}$ and
$\N_{*}$ into what we will call \emph{abstract nerve functors}, i.e.\
functors which, in addition to the properties that wee mentioned in
Pr.~\ref{prop:SmpCtAdj} and Pr.~\ref{prop:PrsStHmt}, have the property
that
\begin{enumerate}
\item \label{AbsNvSmi} in the categories involved all \emph{homotopy
    equivalences} are \emph{weak equivalences}, and
\item \label{AbsNvSmii} the functors between them are \emph{homotopy
    equivalences}.
\end{enumerate}

We do this in two different ways.  In the first we endow the
categories $\RelkCat$, $\simp^{k}\cat S$ and $\simp^{k}\Cathat$ ($k\ge
1$) with \emph{Reedy} structures.  In the second we endow them with
\emph{Rezk} structures which have more weak equivalences and which
turn these categories into models for the theory of
$(\infty,k)$-categories.

%--------------------------------------------------------------------
%--------------------------------------------------------------------
\boldcentered{Abstract nerve functors}

%--------------------------------------------------------------------
\subsection{Definition}
\label{def:AbsNvFn}

A functor $f\colon \cat C \to \cat D$ between saturated relative
categories (Df.~\ref{def:saturated}) will be called an
\textbf{abstract nerve functor} if it has the following four
properties:
\begin{enumerate}
\item \label{AbsNvFni} The relative categories $\cat C$ and $\cat D$
  come with \emph{strict homotopies} for which the associated
  \emph{homotopy equivalences} are \emph{weak equivalences}.
\item \label{AbsNvFnii} The functor $f$ is a relative functor which is
  a \emph{homotopy equivalence}.
\item \label{AbsNvFniii} The functor $f$ sends strictly homotopic maps
  to strictly homotopic maps and hence sends \emph{homotopy
    equivalences} to \emph{homotopy equivalences}.
\item \label{AbsNvFniv} The functor $f$ has a \emph{left adjoint} which
  is a \emph{left inverse}.
\end{enumerate}

%--------------------------------------------------------------------
%--------------------------------------------------------------------
\boldcentered{The Reedy structures}

%--------------------------------------------------------------------
\subsection{Definition}
\label{def:ReeStrct}

We endow
\begin{enumerate}
\item \label{ReeStrcti} the category $\cat S$ with the usual model
  structure,
\item \label{ReeStrctii} the category $\Cathat$ with the \emph{Quillen
  equivalent} Thomason structure \cite{T2},
\item \label{ReeStrctiii} the categories $\simp^{k}\cat S$ and
  $\simp^{k}\Cathat$ ($k\ge 1$) with the resulting Reedy model
  structures,
\item \label{ReeStrctiv} the category $\RelCat$ with the \emph{Quillen
  equivalent} Reedy model structure of \cite{BK1}*{6.1} which was
lifted from the Reedy structure on $\simp\cat S$, and
\item \label{ReeStrctv} the categories $\RelkCat$ ($k>1$) with the
  \emph{homotopy equivalent} relative Reedy structures of
  \cite{BK2}*{3.3} which were lifted from the Reedy structures on
  $\simp^{k}\cat S$.
\end{enumerate}

%--------------------------------------------------------------------
\subsection{Proposition}
\label{prop:AbsNvFncs}

\begin{em}
  If the categories $s^{k}\cat S$, $\RelkCat$ and
  $\simp^{k}\Cathat$ ($k\ge 0$) are endowed with the relative
  structures of \ref{def:ReeStrct}, then the functors $\N_{*}$,
  $\simp^{k}\N$ and $w_{*}$ are abstract nerve functors.
\end{em}

\begin{proof}
  In view of Pr.~\ref{prop:SmpCtAdj} and \ref{prop:PrsStHmt} it
  suffices to show that
  \begin{enumerate}
  \item the functors $N_{*}$, $\simp^{k}\N$ ($k\ge 1$) and $w_{*}$ are
    homotopy equivalences, and
  \item the homotopy equivalences in $\simp^{k}\cat S$, $\RelkCat$ and
    $\simp^{k}\Cat$ ($k\ge 0$) are weak equivalences.
  \end{enumerate}

  To prove (i) we note:
  \begin{lettered}
  \item As Dana Latch \cite{La} has shown that the functor $\N$ is a
    homotopy equivalence it readily follows that so are the functors
    $\N_{*}$.
  \item That the functors $\simp^{k}\N$ are homotopy equivalences was
    shown in \cite{BK2}*{3.4}.
  \item That the functors $w_{*}$ are homotopy equivalences then
    follows a), b) and Pr.~\ref{prop:RlCtsmpS} and the observation
    that homotopy equivalences have the two out of three property.
  \end{lettered}
  To prove (ii) we note
  \begin{resumelettered}{4}
  \item As the homotopy equivalences in $\cat S$ are weak
    equivalences, it readily follows that so are the homotopy
    equivalences in $\simp^{k}\cat S$.
  \item That the homotopy equivalences in $\RelkCat$ and
    $\simp^{k}\Cathat$ are also weak equivalences then follows from d)
    and the fact that the functors $\simp^{k}\N$ and $\N_{*}$ preserve
    homotopy equivalences and reflect weak equivalences.\qedhere
  \end{resumelettered}
\end{proof}

%--------------------------------------------------------------------
%--------------------------------------------------------------------
\boldcentered{The Rezk structures}

%--------------------------------------------------------------------
\subsection{Definition}
\label{def:RzkStrct}

In \cite{R1} Charles Rezk constructed a left Bousfield localization of
the Reedy structure on $\simp\cat S$ and showed it to be a model for
the theory of $(\infty,1)$-categories.

Furthermore it was noted in \cite{B} (and a proof thereof can be found
in \cite{Lu}*{\S1}) that iteration of Rezk's construction yields for
every integer $k>1$, a left Bousfield localization of the Reedy
structure on $\simp^{k}\cat S$ which is a model for the theory of
$(\infty,k)$-categories.

We therefore will denote
\begin{enumerate}
\item \label{RzkStrcti} by $\Le\simp^{k}\cat S$ ($k\ge 1$) the
  category $\simp^{k}\cat S$ with this (iterated) Rezk structure,
\item \label{RzkStrctii} by $\Le\simp^{k}\Cathat$ ($k\ge 1$) the
  induced \cite{H}*{3.3.20} \emph{Quillen equivalent} Rezk model
  structure, and
\item \label{RzkStrctiii} by $\Le\RelkCat$ the \emph{Quillen} or
  \emph{homotopy equivalent} Rezk structure lifted \cite{BK2}*{4.2}
  from the Rezk structure on $\simp^{k}\cat S$ (or the Quillen
  equivalent Rezk structure on $\simp^{k}\Cathat$) which categories
  therefore are all models for the theory of $(\infty,k)$-categories.
\end{enumerate}

%--------------------------------------------------------------------
\subsection{Proposition}
\label{prop:AllAbsnv}

\begin{em}
  If the categories $\simp^{k}\cat S$, $\RelkCat$ and
  $\simp^{k}\Cathat$ are endowed with the Rezk structures of
  \ref{def:RzkStrct}, then the functors $\simp^{k}\N$, $w_{*}$ and
  $\N_{*}$ are abstract nerve functors.
\end{em}

\begin{proof}
  This follows from Pr.~\ref{prop:AbsNvFncs} and the fact that the
  Rezk structures have more weak equivalences than the Reedy ones.
\end{proof}

%--------------------------------------------------------------------
%--------------------------------------------------------------------
%--------------------------------------------------------------------
\part{Homotopy pullback and potential homotopy pullback}

\section{Homotopy pullbacks}
\label{sec:HmtpyPlbk}

%--------------------------------------------------------------------
\subsection{Summary}
\label{sec:HmtpyPlbkSumm}

As we are concerned not only with homotopy pullbacks in the
\emph{model categories} $\RelCat$ and $\simp^{k}\Cathat$ ($k\ge 0$),
but  also in the \emph{saturated relative categories} $\RelkCat$
($k>1$) on which we do \emph{not} have a model structure, we will
define homotopy pullback in a more general fashion than is usually
done.
\begin{enumerate}
\item \label{HmtpyPlbkSummi} In a model category we define a homotopy
  pullback of a zigzag as \emph{any object} which is weakly equivalent
  to its image under a \emph{``homotopically correct'' homotopy limit functor}.
\item \label{HmtpyPlbkSummii} In a saturated relative category we then
  define a homotopy pullback of a zigzag as \emph{any object} weakly
  equivalent to its image under what we will call a \emph{weak
    homotopy limit functor} which is a functor which has only some of
  the properties of the above (i) homotopy limit functors.
\item \label{HmtpyPlbkSummiii} Our main result then is a \emph{global
    equivalence lemma} which states that, if $f\colon \cat C \to \cat
  D$ is a homotopy equivalence between saturated relative categories,
  then $\cat C$ \emph{admits weak homotopy limit functors iff $\cat D$
    does}, and in that case $f$ \emph{preserves homotopy pullbacks}.
\end{enumerate}

In view of Df.~\ref{AbsNvFnii} and Pr.~\ref{prop:AbsNvFncs} this
result not only takes care of the notion of homotopy pullback in the
categories $\RelkCat$ ($k>1$), but it enables us, in the proof of our
main result in section \ref{sec:PrfThm}, to lift our results from the
model categories $\simp^{k}\Cathat$ ($k>1$) to the relative categories
$\RelkCat$ ($k>1$).

%--------------------------------------------------------------------
\subsection{Remark}
\label{rem:Arbtr}

The results of this section actually hold for homotopy limit functors
(and dually homotopy colimit functors) on \emph{arbitrary} diagram
categories.

%--------------------------------------------------------------------
\subsection{Remark}
\label{rem:MdlLim}

As our definition of homotopy pullbacks is much less rigid that the
usual ones, it might be more correct to refer to them as \emph{models
  for homotopy limits}.

%--------------------------------------------------------------------
%--------------------------------------------------------------------
\boldcentered{Homotopy pullbacks in model categories\\
  and their left Bousfield localizations}

%--------------------------------------------------------------------
\subsection{Definition}
\label{def:HomELim}

Let $\cat E$ denote the $2$-arrow category $\cdot \to \cdot \gets
\cdot$.

Given a model category $\cat M$, we then mean by a \textbf{homotopy
  $\cat E$-limit functor} on $\cat M$ a ``homotopically correct''
homotopy limit functor
\begin{displaymath}
  \holim^{\cat E}\colon \cat M^{\cat E} \longrightarrow \cat M \Comma
\end{displaymath}
i.e.\ a functor which, as for instance in \cite{DHKS}*{20.1}, sends
every object of $\cat M^{\cat E}$ to a \emph{fibrant} object of $\cat
M$ and every (objectwise) weak equivalence in $\cat M^{\cat E}$ to a
\emph{weak equivalence} in $\cat M$.

It has the following property.
%--------------------------------------------------------------------
\subsection{Proposition}
\label{prop:HolimRtAdj}

\begin{em}
  The functor
  \begin{displaymath}
    \Ho\holim^{\cat E}\colon \Ho(\cat M^{\cat E})
    \longrightarrow \Ho\cat M
  \end{displaymath}
  is a right adjoint of the constant diagram functor $\Ho\cat M \to
  \Ho(\cat M^{\cat E})$.
\end{em}

\begin{proof}
  This is a special case of \cite{DHKS}*{20.2}.
\end{proof}

%--------------------------------------------------------------------
\subsection{Definition}
\label{def:HmtpyPlbk}

Given an object $\cat B \in \cat M^{\cat E}$, we will say that an
object $\cat U \in \cat M$ is a \textbf{homotopy pullback} of $\cat
B$, if $\cat U$ is weakly equivalent to $\holim^{\cat E}\cat B$.

%--------------------------------------------------------------------
%--------------------------------------------------------------------
\boldcentered{Quasi-fibrant objects in left Bousfield localizations\\ of
    left proper model categories}

%--------------------------------------------------------------------
\subsection{Definition}
\label{def:QFib}

Let $\cat M$ be a model category and let $\Le\cat M$ be a left
Bousfield localization of $\cat M$, i.e.\ \cite{H}*{3.3.3} a model
category with the \emph{same} cofibrations but \emph{more} weak
equivalences.

If $\cat M$ is left proper, then an object $\cat D \in \Le\cat M$ will
be called \textbf{quasi-fibrant} if
\begin{enumerate}
\item \label{QFibi} $\cat X$ is weakly equivalent in $\cat M$ to a
  fibrant object in $\Le\cat M$
\end{enumerate}
or equivalently \cite{H}*{3.4.6(1)}
\begin{resumeenumerate}{2}
\item \label{Qfibii} one (and hence every) fibrant approximation of
  $\cat X$ in $\cat M$ is fibrant in $\Le\cat M$.
\end{resumeenumerate}

%--------------------------------------------------------------------
\subsection{Proposition}
\label{prop:QFpb}

Let $\Le\cat M$ be a left Bousfield localization of a left proper
model category $\cat M$.  Then, for every zigzag between quasi-fibrant
objects, its homotopy pullbacks in $\cat M$ are quasi-fibrant in
$\Le\cat M$ and also homotopy pullbacks of this zigzag in $\Le\cat M$.

\begin{proof}
  This follows readily from \cite{H}*{3.4.6(1)} and
  \cite{H}*{19.6.5}.
\end{proof}

%--------------------------------------------------------------------
\boldcentered{Homotopy pullbacks in saturated relative categories}

%--------------------------------------------------------------------
\subsection{Definition}
\label{def:WHLF}

Give a saturated relative category $\cat R$ (Df.~\ref{saturi}), a
\textbf{weak homotopy $\cat E$-limit functor} will be a relative functor
\begin{displaymath}
  \wholim^{\cat E}\colon \cat R^{\cat E} \longrightarrow \cat R
\end{displaymath}
for which the induced functor
\begin{displaymath}
  \Ho\wholim^{\cat E}\colon \Ho(\cat R^{\cat E}) \longrightarrow
  \Ho\cat R
\end{displaymath}
is a right adjoint of the constant diagram functor $\Ho c\colon
\Ho\cat R \to \Ho(\cat R^{\cat E})$.

%--------------------------------------------------------------------
\subsection{Definition}
\label{def:HmptPlbk}

Given an object $\cat B \in \cat R^{\cat E}$, we will say that an
object $\cat U \in \cat R$ is a \textbf{homotopy pullback} of $\cat B$
if $\cat U$ is weakly equivalent to $\wholim^{\cat E}\cat B$.

%--------------------------------------------------------------------
\subsection{Remark}
\label{rem:LocEq}

While for a model category any two homotopy limit functors
(Df.~\ref{def:HomELim}) are naturally weakly equivalent, this need not
be the case for these \emph{weak} homotopy limit functors.  However
they still have the following \emph{local equivalence} property which,
because of our \emph{non-functorial} definition of homotopy pullbacks,
is all we will need.

%--------------------------------------------------------------------
\subsection{Proposition}
\label{prop:HPBInd}

\begin{em}
The notion of homotopy pullbacks does not depend on the choice of weak
homotopy $\cat E$-limit functor as
  \begin{enumerate}
  \item \label{HPBindi}
    for any two such functors
    \begin{displaymath}
      \wholim^{\cat E}_{1} \quad\text{and}\quad
      \wholim^{\cat E}_{2}\colon \cat R^{\cat E}\longrightarrow
      \cat R
    \end{displaymath}
    the induced functors
    \begin{displaymath}
      \Ho\wholim^{\cat E}_{1} \quad\text{and}\quad
      \Ho\wholim^{\cat E}_{2}\colon \Ho(\cat R^{\cat E})
      \longrightarrow \Ho\cat R
    \end{displaymath}
    are naturally isomorphic
  \end{enumerate}
  which implies that
  \begin{resumeenumerate}{2}
  \item \label{HPBindii} for every object $\cat B \in \cat R^{\cat
      E}$, the objects $\wholim^{\cat E}_{1}\cat B$ and $\wholim^{\cat
      E}_{2}\cat B$ are weakly equivalent.
  \end{resumeenumerate}
\end{em}

\begin{proof}
  This follows readily from the uniqueness of adjoints and the
  saturation of $\cat R$.
\end{proof}

\intro
We end with another useful property.
%--------------------------------------------------------------------
\boldcentered{A global equivalence lemma}

%--------------------------------------------------------------------
\subsection{Lemma}
\label{lem:WHLEx}

\begin{em}
  Let $f\colon \cat R_{1} \to \cat R_{2}$ be a homotopy equivalence
  (Df.~\ref{HomReliii}) between saturated relative categories.  Then
  \begin{enumerate}
  \item \label{WHLExi} there exist weak homotopy $\cat E$-limit
    functors on $\cat R_{1}$ iff they exist on $\cat R_{2}$
  \end{enumerate}
  in which case
  \begin{resumeenumerate}{2}
  \item \label{WHLExii} an object $\cat U \in \cat R_{1}$ is a
    homotopy pullback of an object $B \in \cat R_{1}^{\cat E}$ iff the
    object $f\cat U \in \cat R_{2}$ is a homotopy pullback of the
    object $f\cat B \in \cat R_{2}^{\cat E}$.
  \end{resumeenumerate}
\end{em}

\begin{proof}
  If $\wholim^{\cat E}_{1}\colon \cat R_{1}^{\cat E} \to \cat R_{1}$
  is a weak homotopy limit functor and $g\colon \cat R_{2} \to \cat
  R_{1}$ is a homotopy inverse (Df.~\ref{HomReliii}) of $f$, then it
  suffices to show that the following composition is also a weak
  homotopy limit functor
  \begin{displaymath}
    \xymatrix@C=4em{
      {\cat R_{2}}
      & {\cat R_{1}} \ar[l]_-{f}
      & {\cat R_{1}^{\cat E}} \ar[l]_-{\wholim_{1}^{\cat E}}
      & {\cat R_{2}^{\cat E}} \ar[l]_-{g^{\cat E}}
    }
  \end{displaymath}
  To do this we successively note the following.
  \begin{enumerate}
  \item The maps
    \begin{displaymath}
      \Ho f\colon \Ho\cat R_{1}\longrightarrow \Ho\cat R_{2}
      \qquad\text{and}\qquad
      \Ho g\colon \Ho\cat R_{2} \longrightarrow \Ho\cat R_{1}
    \end{displaymath}
    are inverse equivalences of categories and hence are both left and
    right adjoint.
  \item This implies the existence of the sequence of adjunctions
    \begin{displaymath}
      \xymatrix@C=4.5em{
        {\Ho\cat R_{2}} \ar@<0.75ex>[r]^-{\Ho g}
        \ar@<0.25ex>[r];[]^-{\Ho f}
        & {\Ho\cat R_{1}} \ar@<0.75ex>[r]^-{\Ho i_{1}}
        \ar@<0.25ex>[r];[]^-{\Ho \wholim^{\cat E}_{1}}
        & {\Ho(\cat R_{1}^{\cat E})} \ar@<0.75ex>[r]^-{\Ho f^{\cat E}}
        \ar@<0.25ex>[r];[]^-{\Ho g^{\cat E}}
        & {\Ho(\cat R_{2}^{\cat E})}
      }
    \end{displaymath}
  \item The composition of the left adjoints equals the composition
    \begin{displaymath}
      \xymatrix@C=3em{
        {\Ho\cat R_{2}} \ar[r]^-{\Ho g}
        & {\Ho \cat R_{1}} \ar[r]^-{\Ho f}
        & {\Ho\cat R_{2}} \ar[r]^-{\Ho i_{2}}
        & {\Ho(\cat R_{2}^{\cat E})}
      }
    \end{displaymath}
  \item As $(\Ho f)(\Ho g)$ is naturally isomorphic to the identity of
    $\Ho\cat R_{2}$, the composition in (iii) is naturally isomorphic to
    \begin{displaymath}
      \Ho i_{2}\colon \Ho\cat R_{2} \longrightarrow \Ho(\cat
      R_{2}^{\cat E})
    \end{displaymath}
    which implies that this map is a left adjoint of the composition
    of the right adjoints in (iii).\qedhere
  \end{enumerate}
\end{proof}

%--------------------------------------------------------------------
%--------------------------------------------------------------------
\section{Potential homotopy pullbacks}
\label{sec:PotHmPb}

%--------------------------------------------------------------------
\subsection{Summary}
\label{sec:PHPBsum}

For every integer $n \ge 1$ we will in a functorial manner embed every
zigzag $\zigzag{f}{\cat X}{\cat Z}{\cat Y}{g}$ in the categories
$\RelkCat$ and $\simp^{k}\Cathat$ ($k\ge 0$) (Df.~\ref{def:kSmpNv} and
\ref{def:lvlnrv}) in a commutative diagram of the form
\begin{displaymath}
  \xymatrix{
    {\pullback{\cat X}{\cat Z}{\cat Y}} \ar[r]^-{k} \ar[d]
    & {\subcommacat{f\cat X}{n}{g\cat Y}} \ar[r]^-{\pi} \ar[d]
    & {\cat Y} \ar[d]^{g}\\
    {\cat X} \ar[r]^-{h}
    & {\subcommacat{f\cat X}{n}{\cat Z}} \ar[r]^-{\pi}
    & {\cat Z}
  }
\end{displaymath}
in which
\begin{enumerate}
\item the squares are pullback squares,
\item $h$ is a weak equivalence, and
\item the object $\subcommacat{f\cat X}{n}{g\cat Y}$ is a
  \emph{potential} homotopy pullback of this zigzag in the sense that,
  under suitable restrictions on the map $f\colon \cat X \to \cat Z$
  (which will be discussed in section \ref{sec:PrpBnCn}), this object
  is indeed a homotopy pullback of this zigzag.
\end{enumerate}

We start with an auxiliary construction.
%--------------------------------------------------------------------
\boldcentered{$n$-arrow path objects in $\RelkCat$ ($k\ge 0$)}

%--------------------------------------------------------------------
\subsection{Definition}
\label{def:ArrPthOb}

Given an integer $n \ge 1$ and an object $\cat Z \in \RelkCat$ ($k\ge
0$) we denote by $\subcommacat{\cat Z}{n}{\cat Z} \in \RelkCat$ the
\textbf{$n$-arrow path object} which has
\begin{enumerate}
\item as objects the $n$-arrow zigzags
  \begin{displaymath}
    \xymatrix{
      {Z_{n}} \ar@{}[r]|{\spacedcdots}
      & {Z_{2}}
      & {Z_{1}} \ar[l] \ar[r]
      & {Z_{0}}
    }\quad\text{in $w\cat Z$,}
  \end{displaymath}
\item as maps in $w\subcommacat{\cat Z}{n}{\cat Z}$ and
  $v_{i}\subcommacat{\cat Z}{n}{\cat Z}$ ($1 \le i \le k$) the
  commutative diagrams of the form
  \begin{displaymath}
    \vcenter{
      \xymatrix{
        {Z_{n}} \ar@{}[r]|{\spacedcdots} \ar[d]
        & {Z_{2}} \ar[d]
        & {Z_{1}} \ar[l] \ar[r] \ar[d]
        & {Z_{0}} \ar[d]\\
        {Z'_{n}} \ar@{}[r]|{\spacedcdots}
        & {Z'_{2}}
        & {Z'_{1}} \ar[l] \ar[r]
        & {Z'_{0}}
      }
    }
    \qquad \text{in $a\cat Z$}
  \end{displaymath}
  in which the vertical maps are in $w\cat Z$ and $v_{i}\cat Z$
  respectively, and
\item as maps in $a\subcommacat{\cat Z}{n}{\cat Z}$ those commutative
  diagrams as above which are finite compositions of maps in the
  $v_{i}\subcommacat{\cat Z}{n}{\cat Z}$ ($1 \le i \le k$).
\end{enumerate}
Furthermore
\begin{resumeenumerate}{4}
\item we denote by
  \begin{displaymath}
    \cat Z \xrightarrow{\enskip \pi_{n}\enskip}
    \subcommacat{\cat Z}{n}{\cat Z}
    \xrightarrow{\enskip\pi_{0}\enskip} \cat Z
    \qquad\text{and}\qquad
    \cat Z \xrightarrow{\enskip j\enskip} \subcommacat{\cat Z}{n}{Z}
  \end{displaymath}
  the restrictions of $\subcommacat{\cat Z}{n}{\cat Z}$ to the first
  and last entries respectively and the map which sends each object of
  $\cat Z$ to the alternating zigzag of its identity maps.
\end{resumeenumerate}

These maps have the following nice properties.
%--------------------------------------------------------------------
\subsection{Proposition}
\label{prop:HEWeEq}
\begin{em}
  \begin{enumerate}
  \item \label{HEWeEqi} $\pi_{n}j = \pi_{0} j = 1$,
  \item \label{HEWeEqii} $j$ is a homotopy equivalence which has
    $\pi_{n}$ and $\pi_{0}$ as homotopy inverses
  \end{enumerate}
  and hence (Df.~\ref{ReeStrctiv} and (v) and
  Pr.~\ref{prop:AbsNvFncs})
  \begin{resumeenumerate}{3}
  \item \label{HEWeEqiii} all three maps are weak equivalences in
    $\RelkCat$.
  \end{resumeenumerate}
\end{em}

\begin{proof}
  This is a straightforward computation.
\end{proof}

%--------------------------------------------------------------------
\boldcentered{$n$-arrow fibers and pullback objects in $\RelkCat$ ($k
  \ge 0$)}

%--------------------------------------------------------------------
\subsection{Definition}
\label{def:NArrFibr}

Given an integer $n \ge 1$ and a zigzag in $\zigzag{f}{\cat X}{\cat
  Z}{\cat Y}{g}$ in $\RelkCat$ ($k \ge 0$)
\begin{enumerate}
\item \label{NArrFibri} we denote by $\subcommacat{f\cat X}{n}{\cat
    Z}$ the \textbf{$n$-arrow fibers object} which is defined by
  \begin{displaymath}
    \subcommacat{f\cat X}{n}{\cat Z}
    = \pullback{\cat X}{\cat Z}{\subcommacat{\cat Z}{n}{\cat Z}}
    = \lim\bigl(\cat X \xrightarrow{f}\cat Z
    \xleftarrow{\pi_{n}} \subcommacat{\cat Z}{n}{\cat Z}\bigr)
  \end{displaymath}
  and note that it comes with a \textbf{projection map}
  \begin{displaymath}
    \pi\colon \subcommacat{f\cat X}{n}{\cat Z} \longrightarrow \cat Z
  \end{displaymath}
  induced by the map $\pi_{0}\colon \subcommacat{\cat X}{n}{\cat Z}
  \to \cat Z$, and
\item \label{NArrFibrii} we denote by $\subcommacat{f\cat X}{n}{g\cat
    Y}$ the \textbf{$n$-arrow pullback object} which is defined by
  \begin{displaymath}
    \subcommacat{f\cat X}{n}{g\cat Y}
    = \pullback{\subcommacat{f\cat X}{n}{\cat Z}}{\cat Z}{\cat Y}
    = \lim\bigl(\subcommacat{f\cat X}{n}{\cat Z}
    \xrightarrow{\pi} \cat Z \xleftarrow{g} \cat Y\bigr)
  \end{displaymath}
  and note that it comes with a \textbf{projection map}
  \begin{displaymath}
    \pi\colon \subcommacat{f\cat X}{n}{g\cat Y} \longrightarrow \cat Y
  \end{displaymath}
  obtained by ``restriction to $\cat Y$''.
\end{enumerate}

%--------------------------------------------------------------------
\boldcentered{$n$-arrow fibers and pullback objects in
  $\simp^{k}\Cathat$ ($k \ge 0$)}

%--------------------------------------------------------------------
\subsection{Definition}
\label{def:NArrPB}

As (Nt.~\ref{rem:zerorel}) $\RelzCat = \Cathat$, the case
$\simp^{0}\Cathat = \Cathat$ has already been taken care of in
Df.~\ref{def:NArrFibr}.

Given an integer $n \ge 1$ and a zigzag $\zigzag{f}{\cat X}{\cat
  Z}{\cat Y}{g}$ in $\simp^{k}\Cathat$ ($k\le 1$) we can therefore
define the \textbf{$n$-arrow fibers and pullback objects}
$\subcommacat{f\cat Z}{n}{\cat Z}$ and $\subcommacat{f\cat X}{n}{g\cat
  Y}$ and the associated \textbf{projection maps}
\begin{displaymath}
  \pi\colon\subcommacat{f\cat X}{n}{\cat Z} \longrightarrow \cat Z
  \qquad\text{and}\qquad
  \pi\colon\subcommacat{f\cat X}{n}{g\cat Y}\longrightarrow \cat Y
\end{displaymath}
by the requirement that, for every $k$-fold dimension $p_{*} =
(p_{k},\ldots,p_{1})$
\begin{align*}
  \bigl(\subcommacat{f\cat X}{n}{\cat Z}
  \xrightarrow{\pi}\cat Z\bigr)_{p_{*}}
  &= \bigl(\subcommacat{f_{p_{*}}\cat X_{p_{*}}}{n}{\cat Z_{p_{*}}}
  \xrightarrow{\pi} \cat Z_{p_{*}}\bigr) \in \Cathat,
  \quad\text{and}\\
  \bigl(\subcommacat{f\cat X}{n}{g\cat Y}
  \xrightarrow{\pi}\cat Y\bigr)_{p_{*}}
  &= \bigl(\subcommacat{f_{p_{*}}X_{p_{*}}}{n}{g_{p_{*}}
    \cat Y_{p_{*}}} \xrightarrow{\pi} \cat Y_{p_{*}} \bigr) \in
  \Cathat
\end{align*}

\intro
The definitions \ref{def:NArrFibr} and \ref{def:NArrPB} are closely
related as follows.
%--------------------------------------------------------------------
\subsection{Proposition}
\label{prop:wpicom}

\begin{em}
  For every integer $n \ge 1$ and zigzag $\zigzag{f}{\cat X}{\cat
    Z}{\cat Y}{g}$ in $\RelkCat$ ($k\ge 1$) (Df.~\ref{def:HighEq})
  \begin{align*}
    w_{*}\bigl(\subcommacat{f\cat X}{n}{\cat Z} \xrightarrow{\pi}\cat
    Z\bigr) &= \bigl(\commacat{w_{*}fw_{*}\cat X}{w_{*}\cat Z}
    \xrightarrow{\pi} w_{*}\cat Z\bigr) \in \simp^{k}\Cathat,
    \quad\text{and}\\
    w_{*}\bigl(\subcommacat{f\cat X}{n}{g\cat Y}\xrightarrow{\pi}\cat
    Y\bigr) &= \bigl(\subcommacat{w_{*}fw_{*}\cat X}{n}{w_{*}gw_{*}Y}
    \xrightarrow{\pi} w_{*}Y\bigr) \in \simp^{k}\Cathat.
  \end{align*}
\end{em}
\begin{proof}
  This follows readily from the observation that, for every $k$-fold
  dimension $p_{*} = (p_{k},\cdots,p_{1})$ (Df.~\ref{def:ArrPthOb})
  \begin{displaymath}
    w_{p_{*}}\bigl(\subcommacat{\cat Z}{n}{\cat Z}
    \xrightarrow{\pi_{0}} \cat Z\bigr) =
    \bigl(\subcommacat{w_{p_{*}}\cat
      Z}{n}{w_{p_{*}}Z}\xrightarrow{\pi_{0}} w_{p_{*}}Z\bigr) \in
    \Cathat. \qedhere
  \end{displaymath}
\end{proof}

\intro
It remains the pull it all together as we promised in
\ref{sec:PHPBsum}.
%--------------------------------------------------------------------
\subsection{Proposition}
\label{prop:zzEmbed}

\begin{em}
  For every integer $n \ge 1$, every zigzag $\zigzag{f}{\cat X}{\cat
    Z}{\cat Y}{g}$ in $\RelkCat$ or $\simp^{k}\Cathat$ ($k\ge 0$) can
  in a functorial manner be embedded in a commutative diagram of the
  form (Df.~\ref{def:NArrFibr} and \ref{def:NArrPB})
  \begin{displaymath}
    \xymatrix{
      {\pullback{\cat X}{\cat Z}{\cat Y}} \ar[r]^-{k} \ar[d]
      & {\subcommacat{f\cat X}{n}{g\cat Y}} \ar[r]^-{\pi} \ar[d]
      & {\cat Y} \ar[d]^{g}\\
      {\cat X} \ar[r]^-{h}
      & {\subcommacat{f\cat X}{n}{\cat Z}} \ar[r]^-{\pi}
      & {\cat Z}
    }
  \end{displaymath}
  in which $h$ and $k$ send each object $X \in \cat X$ and $(X,Y) \in
  \pullback{\cat X}{\cat Z}{\cat Y}$ to the alternating zigzag of
  identity maps of $fX$, and
  \begin{enumerate}
  \item the squares are pullback squares, and
  \item the map $h$ is a weak equivalence.
  \end{enumerate}
\end{em}

\begin{proof}
  That the square on the right is a pullback square follows from
  Df.~\ref{def:NArrFibr} and \ref{def:NArrPB} and that the one on the
  left is so is a simple calculation.

  That $h$ is a weak equivalence follows readily from
  Pr.~\ref{prop:HEWeEq}.
\end{proof}

%--------------------------------------------------------------------
\section{Properties \texorpdfstring{$B_{n}$}{Bn} and \texorpdfstring{$C_{n}$}{Cn}}
\label{sec:PrpBnCn}

%--------------------------------------------------------------------
\subsection{Summary}
\label{sec:BnCnsum}

In final preparation for the formulation of our main results (in
section \ref{sec:MainRlt}) we recall from \cite{DKS}*{\S6} the notions
of \emph{properties} $B_{n}$ and $C_{n}$ in $\Cathat$ and then extend
these notions to the categories $\RelkCat$ and $\simp^{k}\Cathat$ for
$k \ge 1$.

As these notions in $\Cathat$ are closely related to the
\emph{Grothendieck construction} we start with a discussion of the
latter.
%--------------------------------------------------------------------
\subsection{Definition}
\label{def:GrConst}

Given an object $\cat D \in \Cathat$ and a \emph{not necessarily
  relative} functor $F\colon \cat D \to \Cathat$
(Df.~\ref{relcatiii}), the \textbf{Grothendieck construction} on $F$
is the object $\Gr F \in \Cathat$ which has
\begin{enumerate}
\item \label{GrConsti} as \emph{objects} the pairs $(D,A)$ of objects
  $D \in \cat D$ and $A \in FD$, and
\item \label{GrCoonstii} as \emph{maps} $(D_{1},A_{1}) \to
  (D_{2},A_{2})$ the pairs $(d,a)$ of maps
  \begin{displaymath}
    d\colon D_{1} \longrightarrow D_{2} \in \cat D
    \qquad\text{and}\qquad
    a\colon (FD_{1})A_{1} \longrightarrow A_{2} \in FD_{2}
  \end{displaymath}
  of which the \emph{compositions} are given by the formula
  \begin{displaymath}
    (d',a')(d,a) = \bigl(d'd, a'\bigl((Fd)a\bigr)\bigr).
  \end{displaymath}
\end{enumerate}

Such a Grothendieck construction $\Gr F$ comes with a
\textbf{projection functor}
\begin{displaymath}
  \pi\colon \Gr F \longrightarrow \cat D\in \Cathat
\end{displaymath}
which sends an object $(D,A)$ (resp.\ a map $(d,a)$) to the object $D$
(resp.\ the map $d$) in $\cat D$.

\intro
The usefulness of Grothendieck constructions is due to the following
property which was noted by Bob Thomason \cite{T1}*{1.2}.
%--------------------------------------------------------------------
\subsection{Proposition}
\label{prop:GrCnHmCl}
\begin{em}
  \begin{enumerate}
  \item The Grothendieck construction is a homotopy colimit functor on
    $\Cathat$
  \end{enumerate}
  and hence
  \begin{resumeenumerate}{2}
  \item it is homotopy invariant in the sense that every natural weak
    equivalence between two functors $F_{1},F_{2}\colon \cat D \to
    \Cathat$ induces a weak equivalence $\Gr F_{1} \to \Gr F_{2}$.
  \end{resumeenumerate}
\end{em}

%--------------------------------------------------------------------
\subsection{Example}
\label{ex:GrCn}

Given a map $f\colon \cat X \to \cat Z \in \Cathat$ and an integer $n
\ge 1$
\begin{enumerate}
\item \label{GrCni} denote, for every object $Z \in \cat Z$ by
  (Df.~\ref{def:NArrFibr})
  \begin{displaymath}
    \subcommacat{f\cat X}{n}{Z} \subset
    \subcommacat{f\cat X}{n}{\cat Z} \in \Cathat
  \end{displaymath}
  the category consisting of the objects and maps which end at $Z$ or
  $1_{Z}$, and
\item \label{GrCnii} denote by
  \begin{displaymath}
    \subcommacat{f\cat X}{n}{-}\colon \cat Z \longrightarrow \Cathat
  \end{displaymath}
  the \emph{not necessarily relative} functor (Df.~\ref{relcatiii})
  which sends each object $Z \in \cat Z$ to $\subcommacat{f\cat
    X}{n}{Z}$ and each map $z\colon Z \to Z' \in \cat Z$ to the
  functor $\subcommacat{f\cat X}{n}{Z} \to \subcommacat{f\cat
    X}{n}{Z'}$ obtained by ``composition with $z$''.
\end{enumerate}

Then one readily verifies that
\begin{resumeenumerate}{3}
\item \label{GrCniii} $\subcommacat{f\cat X}{n}{\cat Z} =
  \Gr\subcommacat{f\cat X}{n}{-}$, and
\item \label{GrCniv} $\bigl(\subcommacat{f\cat X}{n}{\cat
    Z}\xrightarrow{\pi} \cat Z\bigr) = \bigl(\Gr\subcommacat{f\cat
    X}{n}{-}\xrightarrow{\pi} \cat Z\bigr)$.
\end{resumeenumerate}

%--------------------------------------------------------------------
\subsection{Example}
\label{ex:GrCnExAn}

Given a map $f\colon \cat X \to \cat Z \in \Cathat$ and an integer $n
\ge 1$
\begin{enumerate}
\item \label{GrCnExAni} denote, for every pair of objects $X \in \cat
  X$ and $Z \in \cat Z$, by (Ex.~\ref{ex:GrCn})
  \begin{displaymath}
    \subcommacat{fX}{n}{Z} = \subcommacat{f\cat X}{n}{Z}
  \end{displaymath}
  the category consisting of the objects and maps which start at $fX$
  or $1_{fX}$, and
\item \label{GrCnExAnii} denote by
  \begin{displaymath}
    \subcommacat{f-}{n}{Z}\colon \cat X\longrightarrow \Cathat
    \qquad\text{or}\qquad
    \subcommacat{f-}{n}{Z}\colon \cat X\op\longrightarrow \Cathat
  \end{displaymath}
  the functor which sends each object $X \in \cat X$ to
  $\subcommacat{fX}{n}{Z}$ and each map $x\colon X \to X' \in \cat X$
  to the induced functor
  \begin{displaymath}
    \subcommacat{fX}{n}{Z}\longrightarrow \subcommacat{fX'}{n}{Z}
    \qquad\text{or}\qquad
    \subcommacat{fX'}{n}{Z}\longrightarrow \subcommacat{fX}{n}{Z}
  \end{displaymath}
  depending on whether $n$ is even or odd.
\end{enumerate}

Then one readily verifies that
\begin{resumeenumerate}{3}
\item \label{GrCnExAniii}
  \begin{displaymath}
    \subcommacat{f\cat X}{n}{Z} = \left\{
    \begin{aligned}
      \Gr\bigl(\subcommacat{f-}{n}{Z}&\colon \cat X
      \longrightarrow \Cathat\bigr)\\
      \text{or}&\\
      \Gr\bigl(\subcommacat{f-}{n}{Z}&\colon \cat X\op
      \longrightarrow \Cathat\bigr)
    \end{aligned}\right.
  \end{displaymath}
\end{resumeenumerate}

%--------------------------------------------------------------------
\boldcentered{Properties $B_{n}$ and $C_{n}$ in $\Cathat$}

%--------------------------------------------------------------------
\subsection{Definition}
\label{def:BnCn}

\begin{enumerate}
\item \label{BnCni} A map $f\colon \cat X \to \cat Z \in \Cathat$ is
  said to have \textbf{property $B_{n}$} ($n\ge 1$) if the functor
  $Ex.~\ref{ex:GrCn}$
  \begin{displaymath}
    \subcommacat{f\cat X}{n}{-}\colon \cat Z \longrightarrow \Cathat
  \end{displaymath}
  is a \emph{relative functor} (Df.~\ref{relcativ}), and
\item \label{BnCnii} an object $\cat Z \in \Cathat$ is said to have
  \textbf{property $C_{n}$} ($n\ge 1$) if every map
  (Df.~\ref{def:smpnv})
  \begin{displaymath}
    \bs 0^{w} \longrightarrow \cat Z \in \Cathat
  \end{displaymath}
  has property $B_{n}$.
\end{enumerate}

\intro
The usefulness of property $C_{n}$ is due to the following result of
\cite{DKS}*{\S6}.
%--------------------------------------------------------------------
\subsection{Proposition}
\label{prop:BnCn}
\begin{em}
  If, given a map $f\colon \cat X \to \cat Z \in \Cathat$, the object
  $\cat Z$ has property $C_{n}$ ($n \ge 1$), then the map $f\colon
  \cat X \to \cat Z$ has property $B_{n}$.
\end{em}

\begin{proof}
  In view of Ex.~\ref{GrCnii} one has to show that every map $z\colon
  Z \to Z' \in \cat Z$ induces a weak equivalence $\subcommacat{f\cat
    Z}{n}{Zj \to \subcommacat{f\cat X}{n}{Z'}}$, or equivalently
  (Ex.~\ref{GrCnExAniii}) a weak equivalence
  \begin{displaymath}
    \Gr\subcommacat{f-}{n}{Z}\longrightarrow
    \Gr\subcommacat{f-}{n}{Z'}
  \end{displaymath}
  But this follows readily from Ex.~\ref{GrCnExAniii} and
  Df.~\ref{GrCnExAnii} and the fact that, in view of property $C_{n}$,
  for every object $X \in \cat X$ the map
  \begin{displaymath}
    \subcommacat{fX}{n}{Z} \longrightarrow \subcommacat{fX}{n}{Z'}
  \end{displaymath}
  is a weak equivalence.
\end{proof}

%--------------------------------------------------------------------
%--------------------------------------------------------------------
\boldcentered{Properties $B_{n}$ and $C_{n}$ in $\simp^{k}\Cathat$ and
$\RelkCat$ ($k\ge 1$)}

%--------------------------------------------------------------------
\subsection{Definition}
\label{def:MpPrpB}

\begin{enumerate}
\item \label{MpPrpBi} A map $f\colon \cat X \to \cat Z \in
  \simp^{k}\Cat$ has \textbf{property $B_{n}$} ($n\ge 1$) if
  \begin{itemize}
  \item for every $k$-fold dimension $p_{*} = (p_{k},\ldots, p_{1})$,
    the map $f_{p_{*}}\colon \cat X_{p_{*}} \to \cat Z_{p_{*}} \in
    \Cathat$ has property $B_{n}$ (Df.~\ref{BnCni}),
  \end{itemize}
\end{enumerate}
and
\begin{resumeenumerate}{2}
\item \label{MpPrpBii} An object $\cat Z \in \simp^{k}\Cat$ has
  \textbf{property $C_{n}$} ($n \ge 1$) if
  \begin{itemize}
  \item for every $k$-fold dimension $p_{*}= (p_{k},\ldots, p_{1})$
    the object $Z_{p_{*}} \in \Cathat$ has property $C_{n}$
    (Df.~\ref{BnCnii}).
  \end{itemize}
\end{resumeenumerate}

%--------------------------------------------------------------------
\subsection{Definition}
\label{def:MpPrBrel}

\begin{enumerate}
\item \label{MpPrBreli} A map $f\colon \cat X \to \cat Z \in \RelkCat$
  has \textbf{property $B_{n}$} ($n\ge 1$) if
  \begin{itemize}
  \item the map $w_{*}f\colon w_{*}\cat X \to w_{*}\cat Z \in
    \simp^{k}\Cathat$ has property $B_{n}$ (Df.~\ref{MpPrpBi}
  \end{itemize}
  or equivalently if
  \begin{itemize}
  \item for every $k$-fold dimension $p_{*} = (p_{k},\ldots, p_{1})$
    the map $w_{p_{*}}f\colon w_{p_{*}}\cat X \to w_{p_{*}}\cat Z
    \in \Cathat$ has property $B_{n}$
  \end{itemize}
\end{enumerate}
and
\begin{resumeenumerate}{2}
\item \label{MpPrBrelii} an object $\cat Z \in \RelkCat$ has
  \textbf{property $C_{n}$} ($n\ge 1$) if
  \begin{itemize}
  \item the object $w_{*}Z \in \simp^{k}\Cathat$ has property $C_{n}$
    (Df.~\ref{MpPrpBii})
  \end{itemize}
  or equivalently if
  \begin{itemize}
  \item for every $k$-fold dimension $p_{*} = (p_{k}, \ldots, p_{1})$
    the object $w_{p_{*}}\cat Z \in \Cathat$ has property $C_{n}$.
  \end{itemize}
\end{resumeenumerate}

%--------------------------------------------------------------------
%--------------------------------------------------------------------
\part{The main results and their proofs}

%--------------------------------------------------------------------
%--------------------------------------------------------------------
\section{The main results}
\label{sec:MainRlt}

%--------------------------------------------------------------------
\subsection{Summary}
\label{sec:MnRsltSum}

Our main results are
\begin{enumerate}
\item \label{MpRsltSumi} \emph{Theorem} \ref{thm:MpPropBn} which
  states that the presence of properties $B_{n}$ and $C_{n}$ ensure
  that the potential homotopy pullbacks of section~\ref{sec:PotHmPb},
  i.e.\ the \emph{$n$-arrow pullback objects}, are indeed
  \emph{homotopy pullbacks}, and
\item \label{MpRsltSumii} \emph{Theorem} \ref{thm:SuffCn} which states
  that the presence of a \emph{strict $3$-arrow calculus} implies
  \emph{property $C_{3}$}.
\end{enumerate}

These results then are applied to $(\infty,1)$-categories and
$(\infty,k)$-categories for $k>1$ to prove
\begin{resumeenumerate}{3}
\item \label{MpRsltSumiii} \emph{Theorem} \ref{thm:QFibQfib} which
  combines Theorem \ref{thm:MpPropBn} and \ref{thm:SuffCn} with
  results from \cite{BK3} to show that \emph{homotopy pullbacks in
    $(\infty,1)$-categories} can be described as \emph{$3$-arrow
    pullback objects} of zigzags between \emph{partial model
    categories}, i.e.\ relative categories which have the \emph{two
    out of six} property and admit a \emph{$3$-arrow calculus}, and
\item \label{MpRsltSumiv} \emph{Theorem} \ref{thm:MpBn} which is an
  application of Theorem \ref{thm:MpPropBn} to $(\infty,k)$-categories
  for $k\ge 1$, which is much weaker than Theorem \ref{thm:QFibQfib},
  because we have no model structure on $\RelkCat$ for $k\ge 1$, nor
  an analog for partial model categories.
\end{resumeenumerate}

We also give proofs of Theorems \ref{thm:SuffCn}, \ref{thm:QFibQfib}
and \ref{thm:MpBn}, but postpone the proof of Theorem
\ref{thm:MpPropBn} until section \ref{sec:PrfThm}.

%--------------------------------------------------------------------
%--------------------------------------------------------------------
\boldcentered{The main results}

%--------------------------------------------------------------------
\subsection{Theorem}
\label{thm:MpPropBn}

\begin{em}
  Given an integer $n \ge 1$, let $\zigzag{f}{\cat X}{\cat Z}{\cat
    Y}{g}$ be a zigzag in $\RelkCat$ or $\simp^{k}\Cathat$ ($k\ge 0$)
  (Df.~\ref{def:kSmpNv} and \ref{def:lvlnrv}) with the property that
  \begin{enumerate}
  \item \label{MpPropBni} the map $f\colon \cat X \to \cat Z$ has
    property $B_{n}$ (Df.~\ref{BnCni}, \ref{MpPrpBi} and
    \ref{MpPrBreli})
  \end{enumerate}
  which (Pr.~\ref{prop:BnCn}) is in particular the case if
  \begin{resumeenumerate}{2}
  \item \label{MpPropBnii} the object $\cat Z$ has property $C_{n}$
    (Df.~\ref{BnCnii}, \ref{MpPrpBii} and \ref{MpPrBrelii}).
  \end{resumeenumerate}
  Then
  \begin{resumeenumerate}{3}
  \item \label{MpPropBniii} the $n$-arrow pullback object
    $\subcommacat{f\cat X}{n}{g\cat Y}$ (Df.~\ref{def:NArrFibr} and
    \ref{def:NArrPB}) is a homotopy pullback (Df.~\ref{def:HmtpyPlbk}
    and \ref{def:HmptPlbk}) of this zigzag.
  \end{resumeenumerate}
\end{em}

The main tools for proving this are the \emph{fibrillations} of
Hopkins and Rezk which we will discuss in section \ref{sec:HopRzF} and
we therefor postpone the proof of Theorem \ref{thm:MpPropBn} until
section \ref{sec:PrfThm}.

%--------------------------------------------------------------------
\subsection{Corollary}
\label{cor:SuffHPB}

\begin{em}
  A sufficient condition in order that the pullback $\pullback{\cat
    X}{\cat Z}{\cat Y}$ of the above zigzag is also a homotopy
  pullback is that the obvious map
  \begin{displaymath}
    \pullback{\cat X}{\cat Z}{\cat Y} \xrightarrow{\enspace k\enspace}
    \subcommacat{f\cat X}{n}{g\cat Y}
  \end{displaymath}
  of Pr.~\ref{prop:zzEmbed} is a weak equivalence.
\end{em}

%--------------------------------------------------------------------
\subsection{Theorem}
\label{thm:SuffCn}

\begin{em}
  A sufficient condition in order that an object $\cat Z \in \RelkCat$
  ($k\ge 0$) has property $C_{3}$ (Df.~\ref{MpPrBrelii}) is that $\cat
  Z$ admits a strict $3$-arrow calculus (Df.~\ref{def:strcalculus} and
  \ref{def:Str3ar}).
\end{em}

\begin{proof}
  \leavevmode
  \begin{enumerate}
  \item The case $k=0$.  This follows from Rk.~\ref{rem:whystrict},
    Df.~\ref{def:strcalculus} and \cite{DK}*{3.3, 6.1, 6.2 and 8.2}.
  \item The case $k\ge 1$.  It follows readily from
    Df.~\ref{def:Str3ar} and \ref{def:HighEq} that, for every $k$-fold
    dimension $p_{*} = (p_{k},\ldots, p_{1})$, the object
    $w_{p_{*}}\cat Z \in \Cathat$ admits a strict $3$-arrow calculus
    and hence (i) has property $C_{3}$.  The desired result then
    follows from Df.~\ref{MpPrBrelii}.\qedhere
  \end{enumerate}
\end{proof}

%--------------------------------------------------------------------
%--------------------------------------------------------------------
\boldcentered{Applications to $(\infty,1)$-categories and
  $(\infty,k)$-categories for $k>1$}

Before formulating the application mentioned in \ref{MpRsltSumiii} we
recall some results from \cite{BK3}.

%--------------------------------------------------------------------
\subsection{Remark}
\label{rem:CSSfib}

Recall from \cite{BK3} that a \textbf{partial model category} is an
object $\cat X \in  \RelCat$ which admits a $3$-arrow calculus
(Df.~\ref{def:calculus}) and has the \textbf{two out of six property}
that, for every three maps $r$, $s$ and $t \in \cat X$ for which the
\emph{two} compositions $sr$ and $ts$ exist and are weak equivalences,
the other $four$ maps $r$, $s$, $t$ and $tsr$ are also weak
equivalences.

It then was shown in \cite{BK3} that
\begin{enumerate}
\item \label{CSSfibi} \emph{for every partial model category $\cat X
    \in \RelCat$, one (and hence every) Reedy fibrant approximation to
  $\simp\N\cat X \in \simp\cat S$ (Df.~\ref{def:smpnv}) is a complete
  Segal space, i.e.\ a fibrant object in $\Le\simp\cat S$
  (Df.~\ref{RzkStrcti})} 
\end{enumerate}
which in view of \cite{BK1}*{6.1} implies that
\begin{resumeenumerate}{2}
\item \label{CSSfibii} \emph{$\cat X$ is a quasi-fibrant object
    (Df.~\ref{def:QFib}) of $\Le\RelCat$ (Df.~\ref{RzkStrctiii})}.
\end{resumeenumerate}
Moreover
\begin{resumeenumerate}{3}
\item \label{CSSfibiii}\emph{every complete Segal space is Reedy
    equivalent to the simplicial nerve of a partial model category.}
\end{resumeenumerate}

%--------------------------------------------------------------------
\subsection{Theorem}
\label{thm:QFibQfib}

If $\zigzag{f}{\cat X}{\cat Z}{\cat Y}{g}$ is a zigzag in $\RelCat$ in
which
\begin{enumerate}
\item $\cat X$ and $\cat Y$ are quasi-fibrant in $\Le\RelCat$ and
  $\cat Z$ is a partial model category,
\end{enumerate}
which in particular is the case if
\begin{resumeenumerate}{2}
\item $\cat X$, $\cat Y$ and $\cat Z$ are all three partial model
  categories,
\end{resumeenumerate}
then
\begin{resumeenumerate}{3}
\item \emph{$\subcommacat{f\cat X}{3}{g\cat Y}$ is a quasi-fibrant
    object of $\Le\RelCat$, which is a homotopy pullback of this
    zigzag not only in $\RelCat$, but also in $\Le\RelCat$.}
\end{resumeenumerate}

\begin{proof}
  This follows readily from Rk.\ref{rem:whystrict},
  Pr.~\ref{prop:QFpb} and Th.~\ref{thm:MpPropBn} and \ref{thm:SuffCn}.
\end{proof}

\intro It remains to deal with the result that was mentioned in
\ref{MpRsltSumiv}.
%--------------------------------------------------------------------
\subsection{Theorem}
\label{thm:MpBn}

\begin{em}
  Let $\zigzag{f}{\cat X}{\cat Z}{\cat Y}{g}$ be a zigzag in
  $\RelkCat$ ($k\ge 1$) for which
  \begin{enumerate}
  \item $\simp^{k}\N\cat X$, $\simp^{k}\N\cat Z$ and $\simp^{k}\N\cat
    Y$ are quasi-fibrant objects of $\Le\simp^{k}\cat S$
    (Df.~\ref{RzkStrcti})
  \end{enumerate}
  or equivalently (Pr.~\ref{prop:AbsNvFncs})
  \begin{resumeenumerate}{2}
  \item $w_{*}\cat X$, $w_{*}\cat Z$ and $w_{*}\cat Y$ are
    quasi-fibrant objects of $\Le\simp^{k}\Cathat$
  \end{resumeenumerate}
  and assume that for some integer $n \ge 1$
  \begin{resumeenumerate}{3}
  \item the map $f\colon \cat X \to \cat Z$ has property $B_{n}$
  \end{resumeenumerate}
  which (Pr.~\ref{prop:BnCn}) in particular is the case if
  \begin{resumeenumerate}{4}
  \item the object $\cat Z$ has property $C_{n}$,
  \end{resumeenumerate}
  then
  \begin{resumeenumerate}{5}
  \item $\subcommacat{f\cat X}{n}{g\cat Y}$ is a homotopy pullback of
    this zigzag not only in $\RelkCat$ but also in $\Le\RelkCat$
    (Df.~\ref{RzkStrctiii}).
  \end{resumeenumerate}
\end{em}

\begin{proof}
  In view of Df.~\ref{def:NArrPB} and \ref{def:MpPrBrel},
  Pr.~\ref{prop:wpicom} and Th.~\ref{thm:MpPropBn}
  $w_{*}\subcommacat{f\cat X}{n}{g\cat Y}$ is a homotopy pullback if
  the zigzag
  \begin{displaymath}
    \zigzag{w_{*}f}{w_{*}\cat X}{w_{*}\cat Z}{w_{*}\cat Y}{w_{*}g}
  \end{displaymath}
  in $\simp^{k}\Cathat$ as well as, in view of Pr.~\ref{prop:QFpb}, in
  $\Le\simp^{k}\Cathat$.

  The desired result now follows readily from the \emph{Global
    equivalence lemma} \ref{lem:WHLEx} and the fact that
  (Df.~\ref{def:AbsNvFn} and Pr.~\ref{prop:AbsNvFncs}) $w_{*}$ is a
  \emph{homotopy equivalence}.
\end{proof}

%--------------------------------------------------------------------
%--------------------------------------------------------------------
\section{Hopkins-Rezk fibrillations}
\label{sec:HopRzF}

%--------------------------------------------------------------------
\subsection{Summary}
\label{sec:HopRzFSum}

Our proof of Theorem~\ref{thm:MpPropBn} (in section \ref{sec:PrfThm})
will consist of two parts.  The first consists of a proof for the
\emph{model} categories $\simp^{k}\Cathat$ and $\RelCat$.  In the
second we lift these results for the model categories
$\simp^{k}\Cathat$ ($k>1$) to the \emph{relative} categories
$\RelkCat$ by means of the \emph{Global equivalence lemma}
\ref{lem:WHLEx}.

Our aim in this section is to describe three lemmas which we will need
for the first part.  They involve, each in a different way, the
\emph{fibrillations} of Hopkins and Rezk as follows.
\begin{enumerate}
\item \label{HopRzFSumi} The first lemma is a \emph{Quillen
    fibrillation lemma} which will produce the fibrillations to get us
  started.

  It is essentially a reformulation in terms of relative functors and
  fibrillations as well as a slight strengthening of the lemma that
  Quillen used to prove his Theorem~B and states that
  \begin{itemize}
  \item given an object $\cat D \in \Cathat$, a functor $f\colon \cat
    D \to \Cathat$ is a \emph{relative} functor iff the projection
    functor from its Grothendieck construction to $\cat D$
    (Df.~\ref{def:GrConst}) is a \emph{fibrillation} in $\Cathat$.
  \end{itemize}
\item \label{HopRzFSumii} The second lemma is a \emph{Fibrillation
    lifting lemma} which enables us to obtain more fibrillation as it
  provides
  \begin{itemize}
  \item a sufficient condition on a relative functor in order that it
    \emph{reflects fibrillations}.
  \end{itemize}
\item \label{HopRzFSumiii} The third lemma is a \emph{Hopkins-Rezk
    fibrillation lemma} which shows how, in a right proper model
  category, some of the fibrillations obtained in (i) and (ii) can be
  used to construct homotopy pullbacks.
\end{enumerate}

%--------------------------------------------------------------------
%--------------------------------------------------------------------
\boldcentered{Fibrillations}

We start with recalling from \cite{R2}*{\S2} the notion of what
Charles Rezk called \emph{sharp maps} but which, because of their
fibration-like properties (see \ref{prop:PBfibril} below), we prefer
to call \emph{fibrillations}.

%--------------------------------------------------------------------
\subsection{Definition}
\label{def:fibril}

Given a relative category $\cat R$ \textbf{with pullbacks} (i.e.\
which is closed under pullbacks) a map $p \in \cat R$ is called a
\textbf{fibrillation} if every diagram in $\cat R$ of the form
\begin{displaymath}
  \xymatrix{
    {\cdot} \ar[r]_{j} \ar[d]
    & {\cdot} \ar[r] \ar[d]
    & {\cdot} \ar[d]^{p}\\
    {\cdot} \ar[r]_{i}
    & {\cdot} \ar[r]
    & {\cdot}
  }
\end{displaymath}
in which
\begin{itemize}
\item the squares are pullback squares, and
\item $i$ is a weak equivalence
\end{itemize}
has the property that
\begin{itemize}
\item $j$ is also a weak equivalence.
\end{itemize}

This definition readily implies that, just like the \emph{fibrations}
in a \emph{right proper} model category, fibrillations have the
following properties.

%--------------------------------------------------------------------
\subsection{Proposition}
\label{prop:PBfibril}

\begin{em}
  \begin{enumerate}
  \item \label{PBfibrili} Every pullback of a fibrillation is again a
    fibrillation, and
  \item \label{PBfibrilii} every pullback of a weak equivalence along
    a fibrillation is again a weak equivalence.
  \end{enumerate}
\end{em}

\begin{proof}
  This is straightforward.
\end{proof}

%--------------------------------------------------------------------
%--------------------------------------------------------------------
\boldcentered{The fibrillation lifting lemma}

%--------------------------------------------------------------------
\subsection{Lemma}
\label{lem:FibLftLm}

\begin{em}
  Let $f\colon \cat R_{1} \to \cat R_{2}$ be a relative functor
  between relative categories with pullbacks which
  \begin{itemize}
  \item preserves pullbacks (e.g.\ is a right adjoint), and
  \item reflects weak equivalences.
  \end{itemize}
  Then
  \begin{itemize}
  \item it also reflects fibrillations.
  \end{itemize}
\end{em}

\begin{proof}
  Given a pullback diagram in $\cat R_{1}$ of the form
  \begin{displaymath}
    \xymatrix{
      {\cdot} \ar[r]_{j} \ar[d]
      & {\cdot} \ar[r] \ar[d]
      & {\cdot} \ar[d]^{f}\\
      {\cdot} \ar[r]_{i}
      & {\cdot} \ar[r]
      & {\cdot}
    }
  \end{displaymath}
  in which $i$ is a weak equivalence and of which the image in $\cat
  R_{2}$ is a similar diagram in which $fj$ is a fibrillation.  Then
  $fj$ is a weak equivalence and hence so is $j$.
\end{proof}

%--------------------------------------------------------------------
\subsection{Example}
\label{ex:FibLm}

Examples of functors which satisfy the conditions of the fibrillation
lifting lemma \ref{lem:FibLftLm} are the abstract nerve functors of
Df.~\ref{def:AbsNvFn} and Pr.~\ref{prop:AbsNvFncs}.

%--------------------------------------------------------------------
%--------------------------------------------------------------------
\boldcentered{The Hopkins-Rezk fibrillation lemma}

We recall from \cite{R2}*{\S2} the following
%--------------------------------------------------------------------
\subsection{Lemma}
\label{lem:HPBzz}

\begin{em}
  Let
  \begin{displaymath}
    \xymatrix{
      & {\cat A} \ar[r] \ar[d]
      & {\cat Y} \ar[d]^{g}\\
      {\cat X} \ar[r]^{h}
      & {\cat B} \ar[r]^{\pi}
      & {\cat Z}
    }
  \end{displaymath}
  be a diagram in a right proper model category in which
  \begin{itemize}
  \item the square is a pullback square
  \item $h$ is a weak equivalence, and
  \item $\pi$ is a fibrillation.
  \end{itemize}
  Then $\cat A$ is a homotopy pullback of the zigzag
  \begin{displaymath}
    \zigzag{\pi}{\cat B}{\cat Z}{\cat Y}{g}
  \end{displaymath}
  and hence (Df.~\ref{def:HmtpyPlbk}) also of the zigzag
  \begin{displaymath}
    \zigzag{\pi h}{\cat X}{\cat Z}{\cat Y}{g}
  \end{displaymath}
\end{em}

%--------------------------------------------------------------------
%--------------------------------------------------------------------
\boldcentered{The Quillen fibrillation lemma}

%--------------------------------------------------------------------
\subsection{Lemma}
\label{lem:QFLm}

\begin{em}
  Given an object $\cat D \in \Cathat$ (Df.~\ref{kzeroi}) and a
  functor $F\colon \cat D \to \Cathat$ (Df.~\ref{relcatiii}), the
  following three statements are equivalent.
  \begin{enumerate}
  \item \label{QFLmi} $F$ is a relative functor (Df.~\ref{relcativ}),
  \item \label{QFLmii} the map $\N\pi\colon \N\Gr F \to \N\cat D \in
    \cat S$ (Nt.~\ref{kzeroiii} and Df.~\ref{def:GrConst}) is a
    fibrillation, and
  \item \label{QFLmiii} the map $\pi\colon \Gr F \to \cat D \in
    \Cathat$ is a fibrillation.
  \end{enumerate}
\end{em}

\begin{proof}
  To prove (i)$\Rightarrow$(ii) we note that
  \begin{itemize}
  \item Quillen's proof of this lemma \cite{Q}*{\S1} implies that, for
    every integer $p \ge 0$ and map $\Delta[p]= \N\cat p \to \N\cat
    D$, the pullback
    \begin{displaymath}
      \pullback{\N\cat p}{\N\cat D}{\N\Gr F}
    \end{displaymath}
    of the zigzag $\N\cat p \to \N\cat D\gets \N\Gr F$ is a homotopy
    pullback
  \end{itemize}
  and that
  \begin{itemize}
  \item in view of \cite{R2}*{4.1(i and (ii))} this implies that the
    map $\N\pi \to \N\Gr F \to \N D \in \cat S$ is a fibrillation.
  \end{itemize}

  That (ii)$\Rightarrow$(iii) then follows from the \emph{Fibrillation
    lifting lemma} \ref{lem:HPBzz} and Ex.~\ref{ex:FibLm}.

  Finally (iii)$\Rightarrow$(i) follows by a simple calculation from
  the fact that, in view of the \emph{Hopkins-Rezk fibrillation lemma}
  \ref{lem:HPBzz} the pullbacks of the form
  (Df.~\ref{def:strcalculus})
  \begin{displaymath}
    \pullback{\bs 0^{w}}{\cat D}{\Gr F}
    \qquad\text{and}\qquad
    \pullback{\bs 1^{w}}{\cat D}{\Gr F}
  \end{displaymath}
  are both homotopy pullbacks.
\end{proof}

%--------------------------------------------------------------------
%--------------------------------------------------------------------
\section{A proof of Theorem \ref{thm:MpPropBn}}
\label{sec:PrfThm}

%--------------------------------------------------------------------
\subsection{Preliminaries}
\label{sec:Prelim}

We start with recalling from Pr.~\ref{prop:zzEmbed} that
\begin{enumerate}
\item \label{Prelimi} every zigzag $\zigzag{f}{\cat X}{\cat Z}{\cat
    Y}{g}$ in $\simp^{k}\Cathat$ ($k \ge 0$) or $\RelCat$ can, for
  every integer $n \ge 1$, in a functorial manner be embedded in a
  commutative diagram of the form (Df.~\ref{def:NArrFibr} and
  \ref{def:NArrPB})
  \begin{displaymath}
    \xymatrix{
      & {\subcommacat{f\cat X}{n}{g\cat Y}} \ar[r]^{\qquad\pi} \ar[d]
      & {Y} \ar[d]^{g}\\
      {X} \ar[r]^{h\qquad}
      & {\subcommacat{f\cat X}{n}{\cat Z}} \ar[r]^{\qquad\pi}
      & {Z}
    }
  \end{displaymath}
  in which
  \begin{itemize}
  \item the square is a pullback square, and
  \item $h$ is a weak equivalence
  \end{itemize}
\end{enumerate}
and note that, in view of the \emph{Hopkins-Rezk fibrillation lemma}
\ref{lem:HPBzz}
\begin{resumeenumerate}{2}
\item \label{Prelimii} if the map $\pi\colon \commacat{f\cat X}{\cat
    Z} \to \cat Z$ is a fibrillation, then the object
  $\subcommacat{f\cat X}{n}{g\cat Y}$ is a homotopy pullback of that
  zigzag.
\end{resumeenumerate}

%--------------------------------------------------------------------
\subsection{A proof for the category $\Cathat$}
\label{sec:PrfCtht}

\begin{enumerate}
\item \label{PrfCthti} It follows from Df.~\ref{def:GrConst} and
  \ref{BnCni}, Ex.~\ref{GrCniii} and the \emph{Quillen fibrillation
    lemma} \ref{lem:QFLm} that the map $\pi\colon \subcommacat{f\cat
    X}{n}{\cat Z} \to \cat Z$ is a fibrillation, and
\item \label{PrfCthtii} the desired result now follows from
  \ref{sec:Prelim}.
\end{enumerate}

%--------------------------------------------------------------------
\subsection{A proof for the categories $\simp^{k}\Cathat$ ($k\ge 1$)}
\label{sec:PrfSmpk}

\begin{enumerate}
\item In view of \ref{PrfCthti} and Df.~\ref{def:MpPrpB}, for every
  $k$-fold dimension $p_{*}$, the map
  \begin{displaymath}
    \pi_{p_{*}}\colon
    \subcommacat{f_{p_{*}}\cat X_{p_{*}}}{n}{\cat Z_{p_{*}}}
    \longrightarrow \cat Z_{p_{*}} \in \Cathat
  \end{displaymath}
  is a fibrillation, and
\item if $\prod_{p_{*}} p_{*}$ denotes the product of these maps for
  all $k$-fold dimensions and $\prod_{p_{*}} \Cathat$ denotes the
  corresponding product of copies of $\Cathat$, then clearly the same
  holds for the map
  \begin{displaymath}
    \prod_{p_{*}}\bigl(\pi_{p_{*}}\colon
        \subcommacat{f_{p_{*}}\cat X_{p_{*}}}{n}{\cat Z_{p_{*}}}
        \to \cat Z_{p_{*}}\bigr) \in \prod_{p_{*}} \Cathat
  \end{displaymath}
\item Moreover one readily verifies that, in view of
  Df.~\ref{def:NArrPB}, the obvious map $\simp^{k}\Cathat \to
  \prod_{p_{*}}\Cathat$ satisfies the conditions of the
  \emph{Fibrillation lifting lemma} \ref{lem:FibLftLm}, which implies
  that the map $\pi\colon \subcommacat{f\cat X}{n}{\cat Z} \to \cat Z
  \in \simp^{k}\Cathat$ is also a fibrillation.
\item The desired result now follows from \ref{sec:Prelim}.
\end{enumerate}

%--------------------------------------------------------------------
\subsection{Two proofs for the category $\RelCat$}
\label{sec:TwoPrf}

\begin{enumerate}
\item \label{TwoPrfi} As the map $f\colon \cat X \to \cat Z \in
  \RelCat$ has property $B_{n}$ (by assumption), so does, in view of
  Df.~\ref{def:MpPrBrel}, the map $w_{*}f\colon w_{*}\cat X \to
  w_{*}\cat Z \in \simp^{k}\Cathat$.
\item \label{TwoPrfii} It follows from \ref{sec:PrfSmpk} that the map
  $\pi\colon \commacat{w_{*}fw_{*}\cat X}{w_{*}\cat Z} \to w_{*}\cat Z
  \in \simp^{k}\Cathat$ is a fibrillation.
\item \label{TwoPrfiii} Moreover as (Ex.~\ref{ex:FibLm}) $w_{*}$
  satisfies the conditions of the \emph{Fibrillation lifting lemma}
  \ref{lem:FibLftLm} it follows from \ref{prop:wpicom} that the map
  $\pi\colon \subcommacat{f\cat X}{n}{\cat Z} \to \cat Z \in \RelCat$
  is also a fibrillation and the desired result then follows from
  \ref{sec:Prelim}.
\end{enumerate}

However if one does not want to use the model structure on $\RelCat$
one can also, instead of using (ii) and (iii) proceed as follows:
\begin{itemize}
\item [(ii)$'$] It follows from \ref{sec:PrfSmpk} that the object
  $\commacat{w_{*}fw_{*}\cat X}{w_{*}gw_{*}\cat Y}$ is a homotopy
  pullback of the zigzag $\zigzag{w_{*}f}{w_{*}\cat X}{w_{*}\cat
    Z}{w_{*}\cat Y}{w_{*}g}$, and
\item [(iii)$'$] as (Pr.~\ref{prop:AbsNvFncs}) $w_{*}$ is a homotopy
  equivalence the desired result now follows from Pr.~\ref{prop:wpicom}
  and the \emph{Global equivalence Lemma} \ref{lem:WHLEx}.
\end{itemize}

%--------------------------------------------------------------------
\subsection{A proof for the categories $\RelkCat$ ($k>1$)}
\label{sec:PrfRlkCt}

This is essentially the same as the second proof for $\RelCat$
(\ref{TwoPrfi}, (ii)$'$ and (iii)$'$).

%--------------------------------------------------------------------
%--------------------------------------------------------------------

\end{document}